%August, 2002
\input amstex
\documentstyle{amsppt}
\topmatter
\magnification=1200 
\hsize=13.8cm
\vsize=560pt 
\def\R{\Bbb R}
\def\D{\Cal D}
\def\p{\partial}
\def\ol{\overline}
\def\lf{\left}
\def\ri{\right}
\def\Si{\Sigma}
\def\Sir{{\Sigma_{\rho}}}
\def\Sio{{\Sigma_0}}
\def\R{\Bbb R}
\def\tn{\widetilde\nabla}
\def\e{\epsilon} 
\def\wt{\widetilde}
\def\D{\Cal D}
\def\p{\partial}
\def\z{\zeta}
\def\ol{\overline}
\def\b{\beta}
\def\ab{{\alpha\beta}}
\def\lf{\left}

\def\ri{\right}

\def\Si{\Sigma}

\def\Sir{{\Sigma_{\rho}}}

\def\Sio{{\Sigma_0}}

\def\tD{\widetilde\Delta}

\def\tm{\widetilde m}

\def\tg{\tilde g}

\def\Rn{\Bbb R^n}

\def\bx{\bold X}

\def\bn{\bold N}

\def\sigr{\Sigma_r}

\def\By{\bold Y}

\def\aint{\frac{\ \ }{\ \ }{\hskip -0.4cm}\int}

\def\aaint{\frac{\ \ }{\ \ }{\hskip -0.35cm}\int}

\def\gr{\ _rg}

\def\tgr{\ _r\tg}

\def\a{\alpha}

\def\tn{\widetilde\nabla}
\leftheadtext{ Yuguang Shi and Luen-fai Tam}
\rightheadtext{Manifolds with boundary and with  nonnegative scalar curvature}

\topmatter

\title{Positive mass theorem and the boundary behaviors of compact manifolds with nonnegative scalar curvature}\endtitle

\author{Yuguang Shi\footnotemark and  Luen-Fai Tam\footnotemark}\endauthor
 \footnotetext"$^1$"{Research partially supported by NSF of China, Project number 10001001.}
\footnotetext"$^2$"{Research partially supported by   Earmarked
Grant of Hong Kong \#CUHK4032/02P.}
\address{Yuguang Shi: Key Laboratory of  Pure and Applied Mathematics, School of Mathematics Science,
Peking University, Beijing, 100871, China}
\endaddress
\email{ygshi\@math.pku.edu.cn}
\endemail
\address
{Luen-fai Tam:Department of Mathematics, The Chinese
University of Hong Kong, Shatin, Hong Kong, China}
\endaddress
\email{lftam\@math.cuhk.edu.hk}
\endemail
\affil {
Peking University\\
The Chinese University of Hong Kong }
\endaffil
\abstract
In this paper, we study the boundary behaviors of compact manifolds with
nonnegative scalar curvature and nonempty boundary. Using a general version of Positive Mass Theorem of Schoen-Yau and Witten, we prove the following theorem:   For any compact manifold with boundary and
nonnegative scalar curvature, if it is spin and its boundary can
be isometrically embedded into Euclidean space as a strictly
convex hypersurface, then the integral of mean curvature of the
boundary of the manifold cannot be greater than  the integral of
mean curvature of the embedded image as a hypersurface in
Euclidean space. Moreover, equality holds if and only if the manifold is
isometric with a domain in the Euclidean space. Conversely, under the assumption that  the theorem is true, then one can prove the ADM mass of an asymptotically flat manifold is  nonnegative, which is part of the Positive Mass Theorem. 
\endabstract
%\date September, 2002\enddate
\endtopmatter

\subheading{\bf \S0 Introduction}

The structure of a manifold with positive or nonnegative scalar curvature has been studied extensively. There are many beautiful results for compact manifolds without boundary, see \cite{L, SY1-2, GW1-3}. For example, in \cite{L}, Lichnerowicz found that some compact
manifolds admit  no Riemannian metric with positive scalar curvature. In \cite{SY1-2} Schoen and Yau proved that  every torus $T^n$ with $n\le 7$ admits no metric with positive scalar curvature, and admits no non-flat metric with nonnegative scalar curvature. This is also proved later by Gromov and Lawson \cite{GW3} for $n>7$.

For complete noncompact manifolds, the most famous result is the Positive Mass Theorem (PMT), first proved by Schoen and Yau \cite{SY3-4} and later by Witten \cite{Wi} using spinors, see also \cite{PT, B1}. One of their results is as follow: Suppose $(M, g)$ is an asymptotically
flat manifold  such that  $g$ behaves  like Euclidean at infinity near each end, and suppose its scalar curvature is nonnegative, then $(M,g)$ is actually flat if the ADM mass of one of the ends is zero.

It is natural to ask what we can say about  manifolds with boundary and with nonnegative scalar curvature. In a recent work of Yau \cite{Y}, it was proved that if $\Omega$ is a noncompact complete three manifold with boundary and with scalar curvature not less than $-3/2c^2$. Suppose one of the component of $\p\Omega$ has nonpositive Euler number and mean curvature is not less than $c$ and suppose $Area(\p B)\ge c \cdot Vol(B)$ for  any ball $B$ in $\Omega$. Then $\Omega$ is a isometric to the warped product of the flat torus with a half line. This is a result on the effect of mean curvature of the boundary that can influence  the internal geometry of a manifold.

In this work, we will  study   boundary behaviors of compact manifolds  with nonnegative scalar curvature. It turns out that the question is related to the Positive Mass Theorem. The results in this work might also be related to the study of   the quasi-local mass   defined in \cite{B2}.

Consider a compact oriented  three manifold $\Omega^3$ with smooth boundary $\p\Omega$. Suppose each component $\Sigma$ of the boundary has positive Gaussian curvature, then $\Sigma$ can be isometrically embedded in $\R^3$. Moreover, the embedding is unique up to an isometry of $\R^3$, see \cite{N, H}, for example. We will prove:

\proclaim{Theorem 1} Let $(\Omega^3,g)$ be a compact manifold of dimension three with smooth boundary and with nonnegative scalar curvature. Suppose $\p\Omega$ has finitely many components $\Sigma_i$ so that each component has positive Gaussian curvature and positive mean curvature $H$ with respect to the unit outward normal.  Then for each boundary component $\Sigma_i$,
$$
\int_{\Sigma_i}Hd\sigma\le \int_{\Sigma_i}H^{(i)}_0d\sigma\tag0.1
$$
where $H^{(i)}_0$ is the mean curvature of $\Sigma_i$ with respect to the outward normal when it is isometrically embedded in $\R^3$, $d\sigma$ is the volume form on $\Sigma_i$ induced from $g$. Moreover, if equality holds in (0.1) for some $\Sigma_i$, then $\p\Omega$ has only one   component and $\Omega$ is a domain in $\R^3$.
\endproclaim

A similar result is still true in higher dimensions if we assume that each component $\Sigma_i$ can be realized as a strictly convex hypersurface in the Euclidean space and if in addition $\Omega$ is spin. See Theorem 4.1 for more details.

The idea of the proof of Theorem 1 is as follows. We use the methods introduced by Bartnik \cite{B3} to glue the manifold $\Omega$ to another one so that the resulting manifold $N$ is asymptotically flat. This can be accomplished as in \cite{B3} (see also \cite{SW}) by solving a parabolic partial differential equation of some foliation, so that the mean curvatures on the boundary of $\Omega$ and $N\setminus \Omega$ match along $\p\Omega$. Note that the manifold $N$ is only Lipschitz. Next, we will prove that the positive mass theorem is still true for such a manifold, see Theorem 3.1. This theorem  is believed to be true, but the authors are unable to find an explicit reference in the literature and it seems the proof involves some technical points.  We will give a detailed proof of the result. After   obtaining $N$, it can be shown that there is a monotonicity on the difference of the integrals of the mean curvatures of the boundary as a submanifold in $\Omega$ and as a submanifold in the Euclidean space. Then one can conclude the theorem is true.

It is interesting to see that in some sense Theorem 1 is equivalent to the positive mass theorem. In fact, we can prove that:

\proclaim{Theorem 2} Suppose (0.1) is true for any compact Riemannian three manifold $\Omega$ with boundary satisfying the assumptions in  Theorem 1.  Let $(N,g)$ be an asymptotically flat manifold (in a certain sense) with nonnegative  scalar curvature which is in $L^1(N)$. Then the ADM mass $m_E$ is nonnegative for each end $E$ of $N$.
\endproclaim

The paper is organized as follows. In \S1, the equation of foliation is derived. In \S2, we will solve the equation of foliation and obtain necessary estimates for later applications. In \S3, we will prove a  of positive mass theorem for a class of manifolds with Lipschitz metrics. Theorem 1 and its higher dimension analog will be proved in \S4. Theorem 2 will be proved in \S5.

The authors would like to thank Professors Robert Bartnik and   Shing-Tung Yau for useful discussions and their  interest in the work. We would also like to thank Professors Hubert L. Bray  Weiyue Ding and Gang Tian for their interest in the work.

\subheading{\bf \S1 The equation of foliation with prescribed scalar curvature}

In this section, we will derive the equation of foliation with prescribed scalar curvature. The equation has been basically obtained in \cite{B3}, see also \cite{SW}. All manifolds in this work are assumed to be orientable.

Let $\Si$  be a smooth compact manifold without boundary with dimension $n-1$ and let $N=[a,\infty)\times\Si$ equipped with a Riemannian metric of the form

$$
ds_0^2=d\rho^2+g_\rho\tag1.1
$$
for a point $(\rho,x)\in N$. Here $g_\rho$ is the induced metric on $\Sir$ which is the level surface $\rho$=constant. Note that for fixed $x\in \Si$, $(\rho,x)$, $a\le \rho<\infty$ is a geodesic. Given a function $\Cal R$ on $N$, we want to find the equation for $u$ such that

$$
ds^2=u^2d\rho^2+g_\rho\tag1.2
$$
has scalar curvature $\Cal R$. Let $\omega_i$, $1\le i\le n-1$ be a local orthonormal coframe on $\Si_0$. Parallel translate $\omega_i$ on  the direction $\frac{\p}{\p\rho}$. Let $\omega_n=d\rho$. Let $e_i$, $1\le i\le n$ be the dual frame of $\omega_i$, and let $\omega_{ij}$ be the connection forms. Then the structure equations of $ds_0^2$ are

$$
d\omega_i=\sum_{j=1}^n\omega_{ij}\wedge\omega_j,\ \
\omega_{ij}+\omega_{ji}=0,
$$
and

$$
d\omega_{ij}-\sum_{k=1}^n\omega_{ik}\wedge\omega_{kj}=-\frac12\sum_{k,l=1}^nR^0_{ijkl}\omega_k\wedge\omega_l.
$$
where $R^0_{ijkl}$ is the curvature tensor with respect to $ds_0^2$. The second fundamental form $h^0_{ij}$, $1\le i,j\le n-1$ of $\Sir$ with respect to the normal $e_n=\frac{\p}{\p \rho}$ is given by

$$
\omega_{ni}=\sum_{j=1}^{n-1}h^0_{ij}\omega_j.\tag1.3
$$

Let $\eta_i=\omega_i$, $1\le i\le n-1$ and let $\eta_n=u\omega_n$.
Then $\eta_i$ is an orthonormal coframe with respect to $ds^2$. Let $\eta_{ij}$ be
the connection forms of $\eta_i$. Direct computations show that
$$
\eta_{ij}=\omega_{ij},\ \ 1\le i,j\le n-1,\tag1.4
$$
and

$$
\eta_{ni}=-(\log u)_i\eta_n+u^{-1}\omega_{ni}, \ \ 1\le i\le
n-1\tag1.5
$$
where $(\log u)_i$ is the derivative of $\log u$ in the $e_i$ direction. In particular, the second fundamental form $h_{ij}$ of $\Sir$ with respect to $ds^2$ is given by

$$
h_{ij}=u^{-1}h^0_{ij}.\tag1.6
$$

We want to compare the curvature tensor $R_{ijkl}$ of $ds^2$ with $R^0_{ijkl}$. For any $1\leq i,j $$\leq n-1$, apply the Gauss equation   to $\Sir$, noting that the metric on $\Sir$ induced by $ds_0^2$ and $ds^2$ are the same, we have:

$$
\split
  R_{ijij}&=R_{ijij}^\rho+h_{ij}^2-h_{ii}h_{jj}\\
&=R_{ijij}^\rho+u^{-2}\lf((h^0_{ij})^2-h_{ii}^0h_{jj}^0\ri)\\
&=R^\rho_{ijij}+u^{-2}\lf(R_{ijij}^0-R_{ijij}^\rho\ri)\\
&=\lf(1-u^{-2}\ri)R_{ijij}^\rho+u^{-2}R_{ijij}^0
\endsplit\tag1.7
$$
where $R_{ijij}^\rho$ is the intrinsic curvature tensor of $\Sir$. To compare $R_{nini}$ with $R^0_{nini}$, we have

$$
\split
&-\frac12 \sum_{k,l=1}^nR_{nikl}\eta_k\wedge \eta_l=d\eta_{ni}-\sum_{k=1}^{n-1} \eta_{nk}\wedge \eta_{ki}\\
&=-\sum_{j=1}^{n-1}(\log u)_{ij}\omega_j\wedge\eta_n-(\log u)_id\eta_n-u^{-2}\sum_{j=1}^n u_j\omega_j\wedge\omega_{ni}\\
&\quad+u^{-1}d\omega_{ni} -\sum_{k=1}^{n-1}\eta_{nk}\wedge\eta_{ki}\\
&=-\sum_{j=1}^{n-1}(\log u)_{ij}\eta_j\wedge\eta_n-\sum_{j=1}^{n-1}(\log u)_i\lf(-(\log u)_j\eta_n+u^{-1}\omega_{nj}\ri)\wedge\eta_j\\&\quad-u^{-2}\sum_{j=1}^n u_j\omega_j\wedge\omega_{ni} +\lf(u^{-1}\sum_{k=1}^{n-1}\omega_{nk}\wedge\omega_{ki} -\sum_{k=1}^{n-1}\eta_{nk}\wedge\eta_{ki}\ri)\\
&\quad-u^{-1}\cdot\frac12\sum_{k,l=1}^nR^0_{nikl}\omega_k\wedge\omega_l\\
&=\text{\rm I+II+III+IV+V} \endsplit\tag1.8
$$
Here $(\log u)_{ij}=e_j\lf(e_i(\log u)\ri)$. Since $\omega_{nj}(e_n)=0$ for all $j$, $\omega_{nj}$ is a linear combination of $\omega_1,\dots,\omega_{n-1}$. The coefficient of $\eta_n\wedge\eta_i$ in II is $ [(\log u)_i]^2.$ By (1.3), the coefficient of $\eta_n\wedge\eta_i$ in III is $-u^{-3}\frac{\p u}{\p\rho} h_{ii}^0.$ Moreover, by (1.5)

$$
u^{-1}\sum_{k=1}^{n-1}\omega_{nk}\wedge\omega_{ki}-\sum_{k=1}^{n-1}\eta_{nk}\wedge\eta_{ki}=\sum_{k=1}^{n-1}
(\log u)_k\eta_n \wedge\eta_{ki}.
$$
The coefficient of $\eta_n\wedge\eta_i$ in IV is:

$$
\sum_{k=1}^{n-1}(\log u)_k\eta_{ki}(e_i).
$$
Hence compare the coefficients of  $\eta_n\wedge\eta_i$ in (1.8), we have

$$
-R_{nini}=(\log u)_{ii}+[(\log u)_i]^2-u^{-3}u_\rho h_{ii}^0+
\sum_{k=1}^{n-1}(\log u)_k\eta_{ki}(e_i) -u^{-2}R^0_{nini}.
$$
Since $$\eta_{ki}(e_i)=\omega_{ki}(e_i)=-<\nabla_{e_i}e_i,e_k>$$,

$$\sum_{i=1}^{n-1}(\log u)_{ii}+
\sum_{k=1}^{n-1}(\log u)_k\eta_{ki}(e_i)=\Delta_\rho \log u.
$$
where $\Delta_\rho$ is the Laplacian on $\Sir$ with respect to the induced metric from $ds_0^2$. Hence

$$
\sum_{i=1}^{n-1}R_{nini}=-u^{-1}\Delta_\rho  u+u^{-3}\frac{\p
u}{\p \rho} H^0+u^{-2}\sum_{i=1}^{n-1}R^0_{nini}.\tag1.9
$$
where $H^0$ is the mean curvature of $\Sir$ with respect to the metric $ds_0^2$. Combining (1.7) and (1.9), the scalar curvature $\Cal R$ of $ds^2$ is given by

$$
\split
&\Cal R=(1-u^{-2})\Cal R^\rho+u^{-2}\sum_{i,j}^{n-1}R^0_{ijij}+2\sum_{i=1}^{n-1}R_{nini}\\
&=(1-u^{-2})\Cal R^\rho+u^{-2}\Cal R^0 -2u^{-1}\Delta_\rho u
+2u^{-3}\frac{\p u}{\p\rho} H^0
\endsplit
$$
where $\Cal R^0$ is the scalar curvature of $N$ with respect to $ds_0^2$ and $\Cal R^\rho$ is the scalar curvature of $\Sir$ with the induced metric. Hence $u^2 d\rho^2+g_\rho$ has   the scalar curvature  $\Cal R$, if and only if $u$ satisfies

$$
H^0\frac{\p u}{\p \rho}=u^2\Delta_\rho u+\frac12(u-u^3)\Cal
R^\rho-\frac12 u\Cal R^0+\frac{u^3}2\Cal R.\tag1.10
$$

{\bf Example 1.} Let $N=\R^3\setminus B(1)$ with the standard Euclidean metric. Then $N=[1,\infty)\times \Si$ where $\Si$ is diffeomorphic to $\Bbb S^2$. The metric on $N$ is given by $d\rho^2+g_\rho$, where $\Si_\rho, g_\rho$ is the standard sphere  with radius $\rho$. Suppose we want to find $u$ with scalar curvature $\Cal R=0$. Then $u$ satisfies:

$$
2\rho^{-1}\frac{\p u}{\p\rho}=u^2\rho^{-2}\Delta_{\Bbb
S^2}u+(u-u^3)\rho^{-2}
$$
where $\Bbb S^2$ is the standard unit sphere. Hence we have

$$
2\rho \frac{\p u}{\p\rho}=u^2 \Delta_{\Bbb S^2}u+(u-u^3).
$$
This is a special form of the equation derived in \cite{B3}.\vskip .1cm

{\bf Example 2.} Let   $\Sio$ be a smooth compact strictly convex
hypersurface in $R^n$. Let $r$ be the distance function from
$\Sio$. Then the metric on the exterior $N$ of $\Sio$ is given by
$dr^2+g_r$, where $g_r$ is the induced metric on $\Si_r$, which is
the hypersurface with distance $r$ from $\Sio$.  The function $u$
with prescribed scalar curvature $\Cal R=0$ is given by

$$
2H_0\frac{\p u}{\p r}=2u^2\Delta_r u+(u-u^3)\Cal R^r
$$
where $H_0$ is the mean curvature of $\Si_r$, $\Cal R^r$ is the scalar curvature of $\Si_r$ and $\Cal R_0$ is the scalar curvature of $\Si_r$ with the induced metric from $\R^n$ and $\Delta_r$ is the Laplacian on $\Si_r$.

{\bf Example 3.} Let $N=\Bbb H^3\setminus B(1)$ with the standard
hyperbolic metric. Then $N=[1,\infty)\times \Si$ where $\Si$ is
diffeomorphic to $\Bbb S^2$. Then the metric on $N$ is given by
$d\rho^2+\sinh^2 \rho g_0$, where $ g_0$ is the standard metric on
the standard unit sphere in ${\Bbb R}^3$. Suppose we want to find
$u$ with scalar curvature $\Cal R=-6$. Then by a direct
computation, we know $u$ satisfies:

$$
\sinh(2\rho) \frac{\p u}{\p\rho}=u^2 \Delta_{\Bbb
S^2}u+(u-u^3)(1+3\sinh^2 \rho).
$$
\subheading{\bf \S2 Solution to the equation of foliation}

In this section, we will solve the equation in Example 2 in \S1. Namely:

Let $\Sio$ be a compact strictly convex hypersurface in $\Rn$, $\bold X$ be the position vector of a point on $\Sio$, and let $\bold N$ be the unit outward normal of $\Sio$ at $\bold X$. Let $\sigr$ be the convex hypersurface described by $\By=\bold X+r\bold N$, with $r\ge0$. The Euclidean space outside $\Sio$ can be represented by
$$
(\Sio\times(0,\infty), dr^2+g_r)
$$
where $g_r$ is the induced metric on $\sigr$. Consider the following initial value problem

$$
\cases 2H_0\frac{\p u}{\p r}&=2u^2\Delta_ru+(u-u^3)\Cal R^r\ \
\text{on $\Sio\times[0,\infty)$}
\\
u(x,0)&=\quad u_0(x)  \endcases \tag2.1
$$
where $u_0(x)>0$ is a smooth function on $\Sigma_0$,  $H_0$ and
$\Cal R^r$ are the mean curvature and scalar curvature of
$\sigr$ respectively, and $\Delta_r$ is the Laplacian operator on
$\sigr$.

We will solve (2.1) and  show that the metric $ds^2=u^2 dr^2+g_r$ is asymptotically flat  outside $\Sio$. We will also compute the mass of $ds^2$. We basically follow the argument in \cite{B3}, see also \cite{SW}. However, some estimates are obtained with different methods.

\proclaim{Lemma 2.1} Let $(x_1,\dots,x_{n-1})$ be  local coordinates on an open set in $\Sio$. For any integer $k \ge0$ and any multi-index $\alpha$ there is a constant $C$ such that for $r\ge1$

$$
\lf|\lf(r\frac{\p}{\p r}\ri)^k \lf(\frac{\p^{|\alpha|}}{\p
x^\alpha}\ri)\lf(H_0(x,r)-\frac{n-1}r\ri)\ri|\le \frac C {r^2}
$$
and
$$
\lf|\lf(r\frac{\p}{\p r}\ri)^k \lf(\frac{\p^{|\alpha|}}{\p
x^\alpha}\ri) \lf(\Cal R^r
(x,r)-\frac{(n-1)(n-2)}{r^2}\ri)\ri|\le \frac{C}{r^3}.
$$
\endproclaim

\demo{Proof} Let $x$ be a point on $\Sio$ and choose  local coordinates
$(x_1,\dots,x_{n-1})$ near $x$ such that $\frac{\p}{\p x_i}$ is orthonormal at $x$ and such that $\frac{\p\bn}{\p x_i}=k_i\frac{\p\bx}{\p x_i}$. Namely, $k_i>0$ are the principal curvatures of $\Sio$ at $\bx$. Direct computations show that at the point $\By=\bx+r\bn$,

$$
H_0-\frac{n-1}r=\frac1r\sum_{i=1}^{n-1}\frac{1}{1+rk_i}=-\frac1r\frac{\sum_{i=0}^{n-2}
b_ir^i}{\sum_{i=0}^{n-1}a_ir^i},\tag2.2
$$
and
$$
\Cal R^r-\frac{(n-1)(n-2)}{r^2}=-\frac1{r^2}\sum_{1\le i,j\le
n-1,\ i\neq
j}\frac{1+rk_i+rk_j}{(1+rk_i)(1+rk_j)}=\frac{1}{r^2}\frac{\sum_{i=0}^{n-2}
d_ir^i}{\sum_{i=0}^{n-1}c_ir^i},\tag2.3
$$
where $a_i,\ b_i,\ c_i,\ d_i$ are smooth functions on $\Sio$, such that $a_{n-1}>0$ and $c_{n-1}>0$.

Now if $(x_1,\dots,x_{n-1})$ are any local coordinates near a point $x_0$, and if

$$
f(x,r)=\frac{\sum_{i=0}^{p} \beta_ir^i}{\sum_{i=0}^{q}\gamma_ir^i}
$$
where $\beta_i$ and $\gamma_i$ are smooth functions on $\Sio$ with $\gamma_q>0$, then for each $j$,

$$
\frac{\p f}{\p x_j}=\frac{\sum_{i=0}^{p+q}
\tilde\beta_ir^i}{\sum_{i=0}^{2q}\tilde\gamma_ir^i}
$$
and
$$
r\frac{\p f}{\p r}=\frac{\sum_{i=0}^{p+q}
\hat\beta_ir^i}{\sum_{i=0}^{2q}\hat\gamma_ir^i}
$$
where $\tilde\beta_i$,  $\tilde\gamma_i$, $\hat\beta_i$ and $\hat\gamma_i$ are smooth functions on $\Sio$ with $\tilde\gamma_{2q}>0$ and $\hat\gamma_{2q}>0$.

Combining these observations with (2.2) and (2.3), the results follow.

\enddemo

Next, we will obtain   preliminary estimates for the upper and lower bounds for the solution $u$ of (2.1).

\proclaim{Lemma 2.2} If $u$ is defined for all $r$, then there is a constant $C$ independent of $r$ such that
$$
|u(x,r)-1|\le   Cr^{2-n}
$$
for $r\ge1$. In fact, if $u$ is defined on $[0,R)$, then for $0\le r<R$, we have
$$
\lf[1+C_2\exp\lf(-\int_0^r\xi(s)ds\ri)\ri]^{-\frac12} \le
u(x,r)\le
\lf[1-C_1\exp\lf(-\int_0^r\varphi(s)ds\ri)\ri]^{-\frac12}
$$
where
$$ \varphi(r)=\min_{x\in\Sio}\frac{\Cal R^r(x,r)}{H_0(x,r)}>0,\hskip 1cm
\psi(r)=\max_{x\in\Sio}\frac{\Cal R^r(x,r)}{H_0(x,r)}>0,$$
$$
C_1=1-\lf(\max_{\Sio}u_0+1\ri)^{-2}, \hskip 1cm
C_2=\lf(\min_{\Sio}u_0\ri)^{-2}-1,
$$
and $\xi(r)=\varphi(r)$ if $\min_{\Sio}u_0\le 1$, $\xi(r)=\psi(r)$ if $\min_{\Sio}u_0> 1$.
\endproclaim 
\demo{Proof}  Let
$$
f(r)=\lf[1-C_1\exp\lf(-\int_0^r\varphi(s)ds\ri)\ri]^{-\frac12}.
$$
Then  $f(0)>u_0(x)$ for all $x\in\Sio$. For any $\lambda>1$, we have
$$
\split
\frac{d}{dr}\lf(\lambda f\ri)&=\frac12\lf(\lambda f-\lambda f^3\ri)\varphi\\
&> \frac 12(\lambda f-\lambda^3f^3)\frac{\Cal R^r}{H_0}
\endsplit
$$
where we have used the fact that $0<C_1<1$ so that $f>1$, the fact
that $\lambda>1$ and the definition of $\varphi$. An application
of the maximum principle then shows that $u\le \lambda f$. Since
$\lambda>1$ is arbitrary, we have $u\le f$. Notice that
$\varphi(r)=(n-2)/r+O(r^{-2})$, it is easy to see $u-1\le C'
r^{2-n}$ for some $C'$, if $u$ is defined for all $r$. 

To obtain the lower bound for $u$. Suppose $\min_{\Sio}u_0\le 1$. Let 
$$
h(r)=\lf[1+C_2\exp\lf(-\int_0^r\varphi(s)ds\ri)\ri]^{-\frac12}.
$$
It is easy to see that $h$ is well-defined, $h(0)=\min_{\Sio}u_0$ and $h<1$. Then
$$
\split
\frac{d h}{dr}&=\frac12 (h-h^3)\varphi \\
&\le \frac12(h-h^3)\frac{\Cal R^r}{H_0}
\endsplit
$$
where we have used the fact that $h-h^3>0$ and the definition of $\varphi$. As before, we can conclude that $h\le u$. 

Suppose $\min_{\Sio}u_0> 1$. Let 
$$g(r)=\lf[1+C_2\exp\lf(-\int_0^r\psi(s)ds\ri)\ri]^{-\frac12}.
$$
Then $g$ is well defined because $-1<C_2<0$. Moreover, $g(0)=\min_{\Sio}u_0$ and $g>1$. 
$$
\split
\frac{d g}{dr}&=\frac12 (g-g^3)\psi \\
&\le \frac12(g-g^3)\frac{\Cal R^r}{H_0}
\endsplit
$$
where we have used the fact that $g-g^3<0$. We can obtain the required lower bound for $u$ as before.

If $u$ is defined for all $r$, we also have   $u-1\ge -C''r^{2-n}$ for some constant $C''>0$ if $r$ is large enough.
\enddemo

Because of Lemma 2.2, we have:
\proclaim{Lemma 2.3} (2.1) has a unique solution $u$ for all $r$
which satisfies the estimates in Lemma 2.2.
\endproclaim

We need some estimates for the metric $g_r$. Basically, we need
the fact that $r^{-2}g_r$ will be asymptotically equal to the
standard metric on $\Bbb S^{n-1}$. Since $\Sio$ is convex, the
Gauss map $\bn:\Sio\to\Bbb S^{n-1}$ is a diffeomorphism. Fix local coordinates $(x_1,\dots,x_{n-1})$ on $\Sio$ so
that $\Sio$ is given by $\bx(x_1,\dots,x_{n-1})$. Then the metric
$_rg=\ _rg_{ij}dx_idx_j$ on $\sigr$ is given by
$$
_rg_{ij}=\ _og_{ij}+r\lf[\langle \bn_i,\bx_j\rangle +\langle
\bn_j,\bx_i\rangle\ri]+r^2b_{ij}
$$
where $ _og_{ij}$ is the metric on $\Sio$ and $b_{ij}$ is the
standard metric on $\Bbb S^{n-1}$ in coordinates
$(x_1,\dots,x_{n-1})$ via the Gauss map $\bn$. Let $
_r\tg_{ij}=r^{-2} \gr_{ij}$. Then we have the following estimate.
The proof is similar to that of  Lemma 2.2.

\proclaim{Lemma 2.4} With the above notations, for any $k\ge0$ and
any multi-index $\alpha$ there is a constant $C$ such that for
$r\ge1$,
$$
\lf|\lf(r\frac{\p}{\p r}\ri)^k \lf(\frac{\p^{|\alpha|}}{\p
x^\alpha}\ri)\lf( _r\tg_{ij}-b_{ij}\ri)(x,r)  \ri|\le \frac C r.
$$
\endproclaim

For $r\ge 1$, let $r=e^t$, and so $\frac{\p}{\p t}=r\frac{\p}{\p r}$. Equation (2.1) becomes
$$
\frac{\p u}{\p t}=(rH_0)^{-1} u^2\tD_r u+\frac12(u-u^3)r\Cal
R^r H_0^{-1}.\tag2.4
$$
where $\tD_r$ is the Laplacian on $\sigr$ with respect to the metric $\ _r\tg_{ij}$.

\proclaim{Lemma 2.5} Let $u$ be the solution of (2.1), then in
local coordinates $(x_1,\dots,x_{n-1})$  on $\Sio$, for any $k$
and $\alpha$, there is a constant $C$ such that
$$
\lf|\lf( \frac{\p}{\p t}\ri)^k \lf(\frac{\p^{|\alpha|}}{\p
x^\alpha}\ri)\lf(u(x,r)-1\ri)\ri|\le   C r^{2-n}.
$$
\endproclaim
\demo{Proof}  In local coordinates
$$
\tD_r=\frac{1}{\sqrt{\tgr}}\frac{\p}{\p
x_i}\lf(\sqrt{\tgr}\tgr^{ij}\frac{\p}{\p x_j}\ri),
$$
where $\tgr=\det(\tgr_{ij})$. Hence
$$
\split
\lf(rH_0\ri)^{-1}u^2\tD_r u&=\frac{\lf(rH_0\ri)^{-1}u^2}{\sqrt{\tgr}}\frac{\p}{\p x_i}\lf(\sqrt{\tgr}\tgr^{ij}\frac{\p u}{\p x_j}\ri)\\
&=\frac{\p}{\p x_i}\lf[\lf(rH_0\ri)^{-1}u^2\tgr^{ij}\frac{\p u}{\p x_j}\ri]-\frac{\p}{\p x_i}\lf[\frac{\lf(rH_0\ri)^{-1}u^2}{\sqrt{\tgr}}\ri]\lf[\sqrt{\tgr}\tgr^{ij}\frac{\p u}{\p x_j}\ri]\\
&=\frac{\p}{\p x_i}\lf[a_i(x,t,\p u)\ri]-a(x,t,\p u)\endsplit
$$
where
$$
a_i(x,t,\vec p)=\lf(rH_0\ri)^{-1}\tgr^{ij}p_j
$$
and
$$
a(x,t,\vec p)=\frac{\p}{\p
x_i}\lf[\frac{\lf(rH_0\ri)^{-1}}{\sqrt{\tgr}}\ri]\sqrt{\tgr}\tgr^{ij}u^2p_j+
2\lf(rH_0\ri)^{-1}u\tgr^{ij}p_ip_j,$$
Here in  $a$, $u$ is considered to be  a given function.
Hence, by Lemmas 2.1, 2.2 and 2.4, for $t=\log r$ large enough,

$$
a_ip_i\ge C|p|^2
$$
$$
|a_i|\le C'|p|
$$
and
$$
|a|\le C''\lf(1+|p|^2\ri) 
$$
for some positive constants $C,\ C',\ C''$ independent of $t$.
By \cite{LSU, Th. V.1.1}, for any $t_0\ge1$, there are constants
$\beta>0$ and $C_1>0$ independent of $t_0$, such that

$$
\frac{\lf|u(x,t)-u(x',t)\ri|}{|x-x'|^\beta}+\frac{\lf|u(x,t)-u(x,t')\ri|}{|t-t'|^{\frac
\beta2}}\le C_1\tag2.5
$$
for all $x\neq x'\in \Sio$ and $t\neq t'$ in $[t_0,t_0+1]$ for some positive constant $C$ independent on $t_0$. Now consider the function $v=u-1$, we have

$$
\frac{\p v}{\p t}-\lf(rH_0\ri)^{-1}u^2\tgr^{ij}\frac{\p^2 v}{\p
x_i\p x_j}+\frac{\lf(rH_0\ri)^{-1}u^2}{\sqrt{\tgr}}\frac{\p}{\p
x_i}\lf(\sqrt{\tgr}\tgr^{ij}\ri)\frac{\p u}{\p x_j}-\frac12
\lf(u^2+u\ri)r\Cal R^r H_0^{-1}v=0.
$$
By Lemmas 2.1, 2.2, 2.4 and (2.5), using the interior Schauder
estimates \cite{LSU, Th. IV.10.1, or Friedmann Th.1, Chap. 4}, the
lemma  follows.
\enddemo

As in \cite{B3}, let

$$m=\frac{1}2r^{n-2}\lf(1-u^{-2}\ri).\tag2.6$$

By Lemma 2.5, it is easy to see that:

\proclaim{Corollary 2.1} With the notations  in Lemma 2.5, in
local coordinates $(x_1,\dots,x_{n-1})$  on $\Sio$, for any $k$
and $\alpha$, there is a constant $C$ such that

$$
\lf|\lf( \frac{\p}{\p t}\ri)^k \lf(\frac{\p^{|\alpha|}}{\p
x^\alpha}\ri)m\ri|\le   C.
$$
\endproclaim

Direct computations show that $m$ satisfies

$$
\frac{\p m}{\p r}=u^2H_0^{-1}\Delta_r m+3
u^4r^{2-n}H_0^{-1}|\nabla_r m|^2+\lf(\frac{n-2}r-\Cal R^r
H_0^{-1}\ri)m \tag2.7
$$
where $\nabla_r$ is the gradient with respect to the metric $\gr_{ij}$. Hence

$$
\frac{\p m}{\p t}=u^2\lf(rH_0\ri)^{-1}\tD_r m+3
u^4r^{-1-n}H_0^{-1}|\tn_r m|^2+\lf(n-2-\Cal R^r rH_0^{-1}\ri)m
\tag2.7'
$$
where $\tn_r$ is the gradient with respect to the metric $\tgr_{ij}$ and $t=\log r$ as before.

Let $\Delta$ and $\nabla$ be the Laplacian and the gradient with respect to  pull back metric on $\Sio$ of the standard metric on $\Bbb S^{n-1}$ through the Gauss map. Let $(x_1,\dots,x_{n-1})$ be local coordinates on $\Sio$ as in the setting of Lemma 2.3.

\proclaim{Lemma 2.6} With the above notations,
$$
\frac{\p m}{\p t}=\frac1{n-1}\Delta m+f(x,t)
$$
where $f(t,x)$ is a function such that in a local coordinates, for any $k$ and $\alpha$, there is a constant $C$ such that
$$
\lf|\lf( \frac{\p}{\p t}\ri)^k \lf(\frac{\p^{|\alpha|}}{\p
x^\alpha}\ri)f(x,t)\ri|\le   Ce^{-t}.\tag2.8
$$
\endproclaim
\demo{Proof} Here and below $f(x,t)$ will denote a function satisfying (2.8), but it may vary from line to line. By Lemma 2.1, 2.2, 2.5, it is easy to see that

$$
u^2(rH_0)^{-1}=\frac1{n-1}+f.\tag2.9
$$
By Lemma 2.4 and Corollary 2.1, we have

$$
\split
\tD_rm&=\tgr^{ij}\frac{\p^2 m}{\p x_i\p x_j}+\frac1{\sqrt{\tgr}}\frac{\p}{\p x_i}\lf(\sqrt{\tgr}\tgr^{ij}\ri)\frac{\p m}{\p x_j}\\
&=\Delta m+\lf(\tgr^{ij}-b^{ij}\ri)\tgr^{ij}\frac{\p^2 m}{\p x_i\p
x_j}\\&\qquad+
\lf[\frac1{\sqrt{\tgr}}\frac{\p}{\p x_i}\lf(\sqrt{\tgr}\tgr^{ij}\ri)-\frac1{\sqrt{b}}\frac{\p}{\p x_i}\lf(\sqrt{b}b^{ij}\ri) \ri]\frac{\p m}{\p x_j}\\
&=\Delta m+f,
\endsplit\tag2.10
$$
where $b=\det(b_{ij})$. Combining (2.9) and (2.10), we have
$$
u^2\lf(rH_0\ri)^{-1}\tD_r m=\frac1{n-1}\Delta m+f.\tag11$$
Similarly, one can prove that

$$
3u^4r^{-1-n}H_0^{-1}|\tn_r m|^2+\lf(n-2-\Cal R^r r
H_0^{-1}\ri)m=f.
$$
By (2.7'), the lemma follows. 
\enddemo

\proclaim{Lemma 2.7} In local coordinates on $\Sio$, there is a constant $m_0$
$$
|m-m_0|+|\nabla m|(x,t)+\lf|\frac{\p m}{\p t}\ri|(x,t)\le Ce^{-t}
$$
for some constant $C$ for all $x,t$.
\endproclaim
\demo{Proof} Let $ a(t)=\aaint_{\Bbb S^{n-1}}m(x,t)$. Here and below, the volume form of $\Bbb S^{n-1}$ is understood  to be the standard one if there is no specification. Let $\tm(x,t)=m(x,t)- a(t)$. Then
$$
\frac{d  a}{dt}=\aint_{\Bbb S^{n-1}}f(x,t) .\tag2.12
$$
where $f$ is the function in Lemma 2.6. In particular, we have

$$
\lf|\frac{d  a}{dt}\ri|\le Ce^{-t}.
$$
Hence

$$
\split
\frac{d}{dt}\int_{\Bbb S^{n-1}} \tm^2&=2\int_{\Bbb S^{n-1}}\tm\frac{\p\tm}{\p t}\\
&=\frac{2}{n-1}\int_{\Bbb S^{n-1}}\Delta \tm +2\int_{\Bbb S^{n-1}}\tm\lf(f-\frac{d   a}{dt}\ri)\\
&\le -\frac{2}{n-1}\int_{\Bbb S^{n-1}}|\nabla \tm|^2+C_1e^{-t}\lf(\int_{\Bbb S^{n-1}}|\tm|^2\ri)^\frac12\\
&\le -2\int_{\Bbb S^{n-1}}|\tm|^2+C_1e^{-t}\lf(\int_{\Bbb
S^{n-1}}|\tm|^2\ri)^\frac12
\endsplit
$$
for some constant $C_1$ independent of $t$ where we have used (2.12),  Lemma 2.6, the fact that $\int_{\Bbb S^{n-1}}\tm=0$ and the first eigenvalue of $\Bbb S^{n-1}$ is $n-1$. From this it is easy to see that

$$
\lf(\int_{\Bbb S^{n-1}}\tm^2\ri)^\frac12\le C_2e^{-t}(t+1)\tag2.13
$$
for some constant $C_2$. On the other hand,

$$
\split
\lf(\frac{\p }{\p t}-\frac1{n-1}\Delta\ri)\tm^2&=2\tm\lf(\frac{\p\tm}{\p t}-\frac{1}{n-1}\Delta \tm\ri)-\frac{2}{n-1}|\nabla \tm|^2\\
&\le 2\tm\lf(f-  \frac{d  a}{dt}\ri)\\
&\le -C_3e^{-t}
\endsplit\tag2.14
$$
for some constant $C_3$ independent of $t$, where we have used Corollary 2.1, Lemma 2.6 and (2.12). Hence we have

$$
\lf(\frac{\p }{\p
t}-\frac1{n-1}\Delta\ri)\lf(\tm^2+C_3e^{-t}\ri)\le0.
$$
Using the mean value equality and (2.13), we have

$$
\tm^2(x,t)\le C_4e^{-t}(t+1).
$$
for some $C_4$ independent of $t$ and $x$. Put this back to (2.14)  and iterate, we conclude that for any $0<\a<1$, there is a constant $C_5$ independent of $x$ and $t$ such that

$$
|\tm|(x,t)\le C_5e^{-\a t}.\tag2.15
$$
Since $\tm$ satisfies:

$$
\frac{\p \tm}{\p t}=\frac1{n-1}\Delta\tm +f-\aint_{\Bbb
S^{n-1}}f\tag2.16
$$
where $f$ is the function in Lemma 2.6, by the interior Schauder estimates \cite{F, Chap 4, Th. 1},  we conclude that for some $\beta>0$
$$
|\tm|_{2+\beta,\Bbb S^{n-1}\times[t,t+1]}\le C_6e^{-\a t}\tag2.17
$$
for some constant $C_6$ independent of $t$. By the definition of $f$ in Lemma 2.6, we have

$$
\frac{\p m}{\p t}=\frac1{n-1}\Delta m +f
$$
where $|f(x,t)|\le Ce^{(-1-\a)t}.$  Hence (2.13) can be improved as

$$
\lf(\int_{\Bbb S^{n-1}}\tm^2\ri)^\frac12\le Ce^{-t}
$$
and (2.14) can be improved as

$$
\lf(\frac{\p }{\p t}-\frac1{n-1}\Delta\ri)\tm^2\le  C e^{-2t}.
$$
Hence, we have

$$
|\tm|(x,t)\le C_7e^{-t} \tag2.18
$$
for some constant independent of $x$ and $t$. Using (2.16), (2.18), Lemma 2.6 and the  interior Schauder estimate, (2.17) can be improved as

$$
|\tm|_{2+\beta,\Bbb S^{n-1}\times[t,t+1]}\le C_8e^{- t}
$$
Use the definition of $\tm$, we conclude that

$$
|\nabla m|(x,t)+\lf|\frac{\p m}{\p t}\ri|\le C_9e^{-t}.
$$
From the fact that $|\frac{da}{dt}|\le Ce^{-t}$, we conclude that
there is a constant $m_0$ such that $|a(t)-m_0|\le Ce^{-t}$.
Combining these with (2.18),  the lemma follows.
\enddemo

\proclaim{Lemma 2.8} Let $z_1,\dots, z_n$ be the standard
coordinates on $\Bbb R^n$ and let
$\rho(z)=\lf(\sum_{i=1}^nz_i^2\ri)^\frac12$. Then

$$
u(z)=1+\frac{m_0}{\rho^{n-2}}+v
$$
where $m_0$ is the constant in Lemma 2.7, and $v$ satisfies:
$$
|v|=O(\rho^{1-n}),
$$
and

$$
|\nabla_0 v|(z)=O\lf(\rho^{-n}(z)\ri),
$$
where $\nabla_0 v$ is the Euclidean gradient of $v$.
\endproclaim

\demo{Proof} It is easy to see that $|v|=O(\rho^{1-n})$ by the
definition of $m$ and $m_0$ and the fact that $|r-\rho|$ is
bounded. Let $\tilde v=u-1+\frac{m_0}{r^{2-n}}$. By Lemma 2.7 and
the definition of $m$ and $\tilde v$, in local coordinates of
$\Sio$, we have

$$
\lf|\frac{\p \tilde v}{\p x_i}\ri|=\lf|\frac{\p u}{\p x_i}\ri|\le
C_1r^{2-n}\lf|\frac{\p m}{\p x_i}\ri|\le C_1r^{1-n}.\tag2.19
$$
Also

$$
\split
r^{n-2}u\frac{\p \tilde v}{\p r}&=r^{n-2}u\frac{\p }{\p r}\lf(u-1-\frac{m_0}{r^{n-2}}\ri)\\
&=u\lf[\frac{\p m}{\p r}-\frac{n-2}2r^{n-3}\lf(1-u^{-2}\ri)+(n-2)r^{-1}m_0\ri]\\
&=u\lf[\frac{\p m}{\p r}-(n-2)r^{-1}(m-m_0)\ri].
\endsplit\tag2.20
$$
By (2.19), (2.20), Lemma 2.7, the fact that $r\sim \rho$ and the fact that $r=e^t$, we have

$$
|\nabla_0 \tilde v|=O(r^{-n}).\tag2.21
$$

If we use the notations in Lemma 2.1, we see that
$$
\nabla_0 r=\bn.
$$

So

$$
\frac{\p r}{\p z_i}=\bn_i=\frac{z_i-x_i}r\tag2.22
$$
where $\bn_i$ is the $i$-th component of $\bn$, and $\bx=(x_1,\dots,x_n)$ is the position vector on $\Sio$. Since

$$ v-\tilde v=\frac{m_0}{\rho^{n-2}}-\frac{m_0}{r^{n-2}}.
$$
Combining (2.21), (2.22)  and the fact that $|r-\rho|$ is bounded, the lemma is proved.
\enddemo

By Lemma 2.4 and Lemma 2.5, we  get $|u-1|=O(r^{2-n})$, $|\nabla_0
u| =O(r^{1-n})$, and $|\nabla^2_0 u|=O(r^{-n})$ by a direct
computation, here, $\nabla_0$ and $\nabla^2_0$ are the gradient and
Hessian operator of the Euclidean metric respectively. If we write

$$
u^2dr^2+g_r=\sum_{i,j} g_{ij}dz_idz_j.
$$
Then direct computations show (see the computations in  (2.24), (2.27) below, for example):

$$
|g_{ij}-\delta_{ij}|+\rho|\nabla_0g_{ij}|+\rho^2|\nabla^2_0
g_{ij}|\le C\rho^{2-n}.\tag2.23
$$
By the result in \cite{B1}, the ADM mass of the metric $ds^2=u^2
dr^2+g_r$ is well defined, because the scalar curvature of $ds^2$
is zero outside a compact set.

\proclaim{Lemma 2.9} The ADM mass of the metric $u^2 dr^2+g_r$ is
equal to $c(n)m_0$, where $c(n)$ is a positive constant depending
on $n$.
\endproclaim
\demo{Proof} Let $z$ be the standard metric on $\Bbb R^n$, and consider the metric

$$
g=u^2dr^2+g_r=dr^2+g_r+(u^2-1)dr^2.
$$
If we write $g=\sum_{i,j}g_{ij}dz_idz_j$, then

$$
g_{ij}=\delta_{ij}+b_{ij}
$$
where $\sum_{i,j}b_{ij}dz_idz_j=(u^2-1)dr^2.$ Hence

$$
b_{ij}=(u^2-1)\frac{\p r}{\p z_i}\frac{\p r}{\p z_j}.\tag2.24
$$
The ADM mass of $g$ is given by

$$
\lim_{\rho\to\infty}\int_{\Bbb S^{n-1}}\lf(\frac{\p g_{ij}}{\p
z_i}-\frac{\p g_{ii}}{\p z_j}\ri)\rho^{n-2} z_j dV_0\tag2.25
$$
where $dV_0$ is the standard metric on $\Bbb S^{n-1}$. By (2.22)

$$
\frac{\p^2 r}{\p z_i\p
z_j}=\frac{\delta_{ij}}r-\frac{z_iz_j}{r^3}+O(r^{-2})=\frac{\delta_{ij}}\rho-
\frac{z_iz_j}{\rho^3}+O(\rho^{-2}).\tag2.26
$$

$$
\frac{\p g_{ij}}{\p z_i}= \frac{\p b_{ij}}{\p z_i}=2u\frac{\p
u}{\p z_i}\frac{\p r}{\p z_i}\frac{\p r}{\p
z_j}+(u^2-1)\lf(\frac{\p^2 r}{\p z_i^2}\frac{\p r}{\p
z_j}+\frac{\p r}{\p z_i}\frac{\p^2 r}{\p z_i\p z_j}\ri).\tag2.27
$$
Now
$$
\split
2u\frac{\p u}{\p z_i}\frac{\p r}{\p z_i}\frac{\p r}{\p z_j}&=2\lf(1+ m_0\rho^{2-n}+O\lf(\rho^{1-n}\ri)\ri)\lf(-(n-2)m_0\rho^{-n}z_i+O\lf(\rho^{-n}\ri)\ri)\\
&\qquad\times\lf(\rho^{-2}z_iz_j+O\lf(\rho^{-2}\ri)\ri)\\
&=2(n-2)m_0\rho^{-n}z_j+O\lf(\rho^{-n}\ri).
\endsplit\tag2.28
$$
Here repeated indices mean summation.

$$
(u^2-1)\lf(\frac{\p^2 r}{\p z_i^2}\frac{\p r}{\p z_j}+\frac{\p
r}{\p z_i}\frac{\p^2 r}{\p z_i\p
z_j}\ri)=2m_0(n-1)\rho^{-n}z_j+O\lf(\rho^{-n}\ri).\tag2.29
$$

$$
\frac{\p g_{ii}}{\p z_j}=\frac{\p h_{ii}}{\p z_j}=2u\frac{\p u}{\p
z_j}\frac{\p r}{\p z_i}\frac{\p r}{\p z_i}+2(u^2-1)\lf(\frac{\p
r}{\p z_i}\frac{\p^2 r}{\p z_i\p z_j}\ri).\tag2.30
$$

$$
\split
2u\frac{\p u}{\p z_j}\frac{\p r}{\p z_i}\frac{\p r}{\p z_i}&=2\lf(1- m_0\rho^{2-n}+O\lf(\rho^{1-n}\ri)\ri)\lf((n-2)m_0\rho^{-n}z_i+O\lf(\rho^{-n}\ri)\ri)\\
&\qquad\times\lf(1+O\lf(\rho^{-2}\ri)\ri)\\
&=2(n-2)m_0\rho^{-n}z_j+O\lf(\rho^{-n}\ri).
\endsplit\tag2.31
$$
$$
2(u^2-1)\lf(\frac{\p r}{\p z_i}\frac{\p^2 r}{\p z_i\p z_j}\ri)
=O\lf(\rho^{-n}\ri).\tag2.32
$$
Combining (2.25), (2.27)-(2.32), the lemma follows.
\enddemo

\proclaim{Lemma 2.10}

$$
\lim_{r\to\infty}\int_{\sigr}\lf(H_0-H\ri)d\sigma_r=\lim_{r\to\infty}\int_{\sigr}H_0(1-u^{-1})d\sigma_r=(n-1)\omega_{n-1}m_0
$$
where $\omega_{n-1}$ is the volume of the standard sphere $\Bbb S^{n-1}$ and $H_0$ and $H$ are the mean curvature of $\sigr$ with respect to the Euclidean metric and the metric $u^2dr^2+g_r$ respectively.
\endproclaim
\demo{Proof} The result follows from (1.6), Lemmas 2.1, 2.8, 2.9  and the definition of $\sigr$.
\enddemo

We can summarize the results in Lemma 2.3, 2.8, 2.9, 2.10 as follows:

\proclaim{Theorem 2.1} The initial value problem (2.1) has a unique solution $u$ on $\Sio\times[0,\infty)$ such that

\roster

\item"{(a)}"

$$
u(z)=1+\frac{m_0}{\rho^{n-2}}+v
$$
where $m_0$ is a constant and $v$ satisfies $|v|=O\lf(\rho^{1-n}\ri)$ and $|\nabla_0 v|=O\lf(\rho^{-n}\ri)$;

\item"{(b)}" The metric $ds^2=u^2 dr^2+g_r$ is asymptotically flat in the sense  of (2.23) with scalar curvature $\Cal R\equiv0$ outside $\Sio$;

\item"{(c)}" The ADM mass $m_{ADM}$ of $ds^2$ is given by

$$
c(n)m_{ADM}=(n-1)\omega_{n-1}m_0=\lim_{r\to\infty}\int_{\sigr}H_0(1-u^{-1})d\sigma_r=\lim_{r\to\infty}\int_{\sigr}\lf(H_0-H\ri)d\sigma_r,
$$
for some positive constant $c(n)$, where $H_0$ and $H$ are the mean curvatures of $\sigma_r$ with respect to the Euclidean metric and $ds^2$ respectively.
\endroster
\endproclaim

 If we let $u_0\equiv k$ for $k\ge 1$, it is easy to see from Lemma 2.2, that the solution $u^{(k)}$ of (2.1) are uniformly bounded on $[a,\infty)$ for all $a>0$. Hence as in \cite{B3}, we can solve (2.1) with initial value $u_0^{-1}=0$. In fact,  by Lemma 2.2, $u_0$ satisfies:

$$
\lf[1-\exp\lf(-\int_0^r\psi(s)ds\ri)\ri]^{-\frac12} \le
u_0(x,r)\le
\lf[1-\exp\lf(-\int_0^r\varphi(s)ds\ri)\ri]^{-\frac12}.
$$
This means that $\Sio$ is a minimal surface with respect to the asymptotically flat metric $u^2dr^2+g_r$. As in \cite{B3}, we have the following:

\proclaim{Lemma 2.11}Let
$$M(r)=\int_{\Sigma_r}H_0 (1-u^{-2})d\sigma_r,$$
then M(r) is nondecreasing
\endproclaim
\demo{Proof} By the Gauss equations, it is easy to see that
$$
\frac{\p H_0}{\p r}=-\sum_{i,j=1}^{n-1}\lf(h_{ij}^0\ri)^2 
$$
where $h_{ij}^0$ is the second fundamental form of $\Sigma_r$ in $\R^n$.
By direct computation, we see:
$$\split
\frac{d}{dr}\int_{\Sigma_r}H_0
(1-u^{-2})d\sigma_r&=\int_{\Sigma_r}({H_0}^2 (1- u^{-2})+ 2 H_0
u^{-3}\frac{\p u}{\p r} + \frac{\p H_0}{\p
r}(1-u^{-2}))d\sigma_r \\
&=2\int_{\Sigma_r}u^{-1}\triangle_{\Sigma_r}u d\sigma_r \\
&=2\int_{\Sigma_r}u^{-2}|\nabla u|^2 d\sigma_r\\
&\geq 0 
\endsplit
$$
where we have used the fact that $u$ satisfies (2.1) and that $H_0^2-\sum_{i,j=1}^{n-1}\lf(h_{ij}^0\ri)^2=\Cal R^r.$

Hence, $M(r)$ is nondecreasing.
\enddemo

Thus, as in \cite{B3, Corollary 1.1} we have:

 \proclaim{Proposition 2.1} Let $u$ be the solution of (2.1)
with initial value $u_0^{-1}=0$. Let $m_{ADM}$ be the ADM mass of
the metric $u^2dr^2+g_r$. Then
$$
m_{ADM}\ge C(n)\int_{\Sio}H_0d\sigma_0
$$
for some positive constant $C(n)$ depending only on $n$.

\endproclaim

\subheading{\bf \S3 Positive mass theorem on manifolds with Lipschitz metric}

In order to prove the main result we need to verify that the positive mass theorem is still true for some manifolds whose metrics may be only Lipschitz. In this section, we always assume that $N^n$ is an orientable  complete noncompact smooth manifold with dimension $n$, such that there is bounded domain $\Omega\subset N$ with smooth boundary $\p\Omega$. We also assume that $N$ is spin (which is always true if $n=3$) and there is a continuous Riemannian metric $g$ on $N$ such that

\roster

\item"{(i)}" $g$ is smooth on $N\setminus \Omega$ and $\ol\Omega$, and is Lipschitz near $\p\Omega$.

\item"{(ii)}" The mean curvatures at $\p\Omega$ with respect to the outward normal and with respect to the metrics $g|_{N\setminus \Omega}$ and $g|_{\ol\Omega}$ are the same.

\item"{(iii)}" $N$ has finitely many ends, each of which is asymptotically Euclidean in the following sense: There is a compact set $K$ containing $\Omega$ such that $N\setminus K=\cup_{i=1}^\ell E_i$. Each $E_i$ is diffeomorphic to $\R^n\setminus B_{R_i}(0)$ and in the standard coordinates in $\R^n$, the metric $g$ satisfies

$$g_{ij}=\delta_{ij}+b_{ij},$$
with

$$\|b_{ij}\|+r\|\p b_{ij}\|+r^2\|\p\p b_{ij}\|=O(r^{{2-n}}) \tag3.1 $$
where $r$ and $\p$ denotes Euclidean distance and the standard gradient operator on $\R^n$, respectively.

\item"{(iv)}" The scalar curvature of $N\setminus \p\Omega$ is nonnegative and is in $L^1(N)$.
\endroster

We should remark that because of (i), the outward unit normal on $\p\Omega$ is well-defined. Moreover, (iii) and (iv) imply that the ADM mass of each end of $N$ is also well-defined by the proof in  \cite{B1}. Explicitly, the ADM mass at each end $E$ is given by

$$
C(n)m_E =\lim_{r\to\infty}\int_{S(r)}\lf(g_{ij,j}-g_{jj,i}\ri)dS_i
$$
where $C(n)$ is a positive constant, $S(r)$ is the Euclidean sphere and $dS_i$ is the normal surface area of $S(r)$.

We have the following:

\proclaim{Theorem 3.1} Let $(N,g)$ as above. Then $m_E\ge0$ for any end $E$ of $N$. Moreover, if the ADM mass of one of ends of $N$ is zero, then $N$ has only one end and $N$ is  flat.

\endproclaim

We will use the argument of Witten \cite{Wi, PT, B1}. Let us first fix some notations. In the following, a local orthonormal frame $e_i,\ 1\le 1\le n$ means that $e_i=a_{ij}\frac{\p}{\p x_j}$ with Lipschitz functions $a_{ij}$ which are smooth on $N\setminus \Omega$ and $\ol\Omega$, where $(x_1,\dots,x_n)$ are  smooth local coordinates. By the assumptions on $g$, we can always find such a local frame near each point.

Let $e_i$ be a local orthonormal frames and $\omega_i$ be the dual 1-forms. Let $\omega_{ij}$ be the connection forms of $g$ and let $\{\sigma_I\}$ be the  orthonormal base of fibers of the spinor bundle $\Cal{S}$ with respect to $\{e_i\}$, $\nabla$ be the covariant derivative on $\Cal{S}$, then we have:

$$\nabla \sigma_I =-\frac{1}{4}\sum_{i,j}\omega_{ij}\otimes e_i\cdot e_j \cdot
\sigma_I,$$
where $"\cdot"$ refers to Clifford multiplication. By the above notations, the Dirac operator can be expressed in the following way:

$$\D=\sum_{i=1}^n e_{i}\cdot \nabla_{e_i}.$$

A spinor $\Psi$ is said to be in $W^{1,2}_{loc}(U)$ in some open set $U$ if near each point $x\in U$  there is a local orthonormal frame $e_i$, such that if $\sigma_I$ is a base for $\Cal S$ as above and

$$\Psi=\sum_{I} \Psi^{I} \sigma_I $$
then $\Psi^I$ is in $W^{1,2}$ near $x$. That is to say, $\Psi^I$ has weak derivatives so that $\Psi^I$ together with its weak derivatives are in $L^2_{loc}$. Note this is well-defined because the transition functions from one orthonormal base to another one are Lipschitz. Note that  it is also meaningful to say that $\Psi$ is locally H\"older or locally Lipschitz. For an open set $U$, the norm $W^{1,2}$ norm $|||\Psi|||$ of a spinor $\Psi$ is defined as

$$
|||\Psi|||^2=\int_{U}||\nabla \Psi||^2+||\Psi||^2.
$$

Near a point $x\in   \p\Omega$, choose an orthonormal  frame $e_i$ such that $e_n=\frac{\p}{\p\rho}$ where $\rho$ is the signed distance function from $\p\Omega$. $\rho>0$ outside $\Omega$ and $\rho<0$ in $\Omega$.   Moreover $e_i$, $1\le i\le n-1$, are chosen so that they are obtained  by parallel translation along the integral curves of $\frac{\p}{\p \rho}$ which are geodesics normal to $\p\Omega$. We call this kind of frame to be an {\it adapted} orthonormal frame. Let $\omega_i$ be the dual of $e_i$ and let   $\omega_{ij}$ be the connection forms. It is easy to see that we have the following:

\proclaim{Lemma 3.1} With the above notations,  $\omega_{ij}(e_k)$ are Lipschitz for $1\le i,j\le n-1$ and for all $k$. Moreover, $\omega_{ij}(e_n)=0$ for all $i, \ j$.

\endproclaim

Under the $adapted$ orthonormal frame, we have:

$$\Cal R=-\frac{\p H}{\p \rho}-\sum_{i,j}{h_{ij}}^2 + {\Cal R}^\rho,
$$
where $h_{ij}$ is components of the second fundamental form and $\Cal R^\rho$ is the scalar curvature of the hypersurface with distance $\rho$ from $\p\Omega$.  Since $H$ matches along $\p\Omega$, it is Lipschitz. By this formula we see that $\Cal R$ is well-defined in the distribution sense. This is important in the proof of the following Lichnerowicz formula.

\proclaim{Lemma 3.2} Let $U$ be a open set of $N$. For any spinor
$\eta\in W^{1,2}_0 (U)$, $\Psi \in W^{1,2}_{loc}(U)$, we have:

$$\int_{U}\langle \D \Psi,\D \eta\rangle
=\int_{U}\langle \nabla \Psi, \nabla \eta \rangle +
\frac{1}{4}\int_{U}\Cal R \langle \Psi,\eta\rangle,$$
where $\Cal R$ is the scalar curvature of $N$.
\endproclaim

\demo{Proof} Let $T=\partial \Omega \cap U$.  Since the metric $g$
is  smooth up to the boundary on $\bar{\Omega}$, by the standard
Lichnerowicz formula applied to   $\Omega \cap U$, we have:

$$\int_{\Omega \cap U}\langle \D \Psi,\D \eta\rangle +\int_{T}\langle  \nu\cdot D \Psi+
\nabla_\nu \Psi, \eta\rangle =\int_{\Omega \cap U}\langle \nabla
\Psi, \nabla \eta \rangle + \frac{1}{4}\int_{\Omega \cap U}\Cal R
\langle \Psi,\eta\rangle, \tag3.2
$$
here  $\nu$ is the outer normal unit vector of $\p\Omega$, see \cite{B1, p.689} for example.

Let $\{e_i\}$ be an adapted frame near a point $x\in T$. Direct computations show that

$$\nu\cdot \D \Psi+
\nabla_\nu \Psi =\nu\cdot \sum_{i=1}^{n-1}e_{i}\cdot
(\nabla_{e_i}\Psi^{I})\sigma_I -\frac{1}{4}\nu\cdot
\sum_{i=1}^{n-1}\sum_{s,t=1}^{n-1}\omega_{st}(e_i)e_{i}\cdot e_s
\cdot e_t \cdot \Psi + \frac{H}{4}\nu\cdot \Psi.
$$
where $H$ is the mean curvature of the level set $\rho$=constant
with respect to $\frac{\p}{\p\rho}$. Hence, we have:

$$\split
&\int_{\Omega \cap U}\langle \D \Psi,\D \eta\rangle
+\int_{T}\nu\cdot \sum_{i=1}^{n-1}e_{i}\cdot
(\nabla_{e_i}\Psi^{I})\sigma_I\\
& -\frac{1}{4}\int_{T}\nu\cdot
\sum_{i=1}^{n-1}\sum_{s,t=1}^{n-1}\omega_{st}(e_i)e_{i}\cdot e_s
\cdot e_t \cdot \Psi  + \int_{T}\frac{H}{4}\nu\cdot \Psi\\
 &=\int_{\Omega \cap U}\langle \nabla \Psi, \nabla \eta \rangle +
\frac{1}{4}\int_{\Omega \cap U}\Cal R \langle \Psi,\eta\rangle.
\endsplit
$$
By the same reason, we have the   formula  on $U\setminus \ol\Omega$,

$$\split
&\int_{U\setminus\ol \Omega}\langle \D \Psi,\D \eta\rangle
-\int_{T}\nu\cdot \sum_{i=1}^{n-1}e_{i}\cdot
(\nabla_{e_i}\Psi^{I})\sigma_I\\
& +\frac{1}{4}\int_{T}\nu\cdot
\sum_{i=1}^{n-1}\sum_{s,t=1}^{n-1}\omega_{st}(e_i)e_{i}\cdot e_s
\cdot e_t \cdot \Psi -\int_{T}\frac{H}{4}\nu\cdot \Psi\\
&=\int_{U\setminus\ol \Omega}\langle \nabla \Psi, \nabla \eta \rangle +
\frac{1}{4}\int_{U\setminus\ol \Omega}\Cal R \langle \Psi,\eta\rangle.
\endsplit
$$
In the above, we have used the fact that the mean curvatures of $T$
in $\Omega$ are equal to that of $N\setminus \Omega$, and the unit outward
normals are in opposite directions. Adding these two equalities,  we see the integrals on $T$ are canceled. Hence, the proof of the lemma is
completed.
\enddemo

Let $\Psi$ be a spinor in $W^{1,2}_{loc}(U)$. $\Psi$ is said to satisfy

$$
\D^2\Psi=0
$$
in the weak sense in an open set $U$ if for any spinor $\Phi\in W^{1,2}_0(U)$,

$$
\int_U\langle \D\Psi,\D\Phi\rangle=0.
$$

Even though $g$ is not smooth, however the coefficients of a weak solution $\Psi$ of $\D^2\Psi=0$ with respect to an adapted frame behave well. Namely, we have:

\proclaim{Lemma 3.3} Suppose $\Psi\in W^{1,2}_{loc}(U)$ satisfies $\D^2\Psi=0$ weakly in an open set $U$. Then $\Phi$ is locally H\"older continuous and $\Psi$ is in $W^{2,2}_{loc}(U)$ in the following sense: (a) if $x\notin \p\Omega$, then $\Psi\in W^{2,2}$ near $x$; (b)   if $x\in \p\Omega$, and if $\{e_i\}$ is an adapted orthonormal frame near $x$ so that $\{\sigma_I\}$ is an orthonormal basis for $\Cal S$ with respect to $\{e_i\}$ and that $\Psi=\sum_I\Psi^I\sigma_I$, then $\Psi^I$ is in $W^{2,2}$ and is H\"older continuous near $x$.
\endproclaim

\demo{Proof}  It is sufficient to study the behavior of $\Psi$
near a point in $\p\Omega$. Let $\Psi=\sum_I\Psi^{I}\sigma_{I}$ as
in case (b). We may assume   $U$ is small enough so that there is an adapted orthonormal frame $e_i$ in $U$.  We claim that   $\Psi^I$ satisfies the following
equations:

$$\triangle \Psi^I +\sum_{J,i} A_{Ji}^{I}e_i\Psi^J
+\sum_JB^{I}_{J}\Psi^{J} =0,
$$
in the weak sense, where, $\|A_{Ji}^{I}\|_{L^{\infty}}
+\|B^I\|_{L^{\infty}} \leq C<\infty$ locally, and $\Delta$ is the
Laplacian for function on $N$. In particular, for each fixed $I$, $\Psi^I$ satisfies the following equation in the weak sense:
$$
\triangle \Psi^I =f
$$
where $f=-\sum_{J,i} A_{Ji}^{I}e_i\Psi^J-\sum_{J} B^{J}_{I}\Psi^{J}.$ Since $\Psi\in W^{1,2}_{loc}(U)$, $f\in L^2_{loc}(U)$. Since the metric is Lipschitz, in local coordinates $\Delta$ is  of divergence form with coefficients being Lipschitz. Then, by the standard theory in
elliptic equations, we know that $\Psi^I\in W^{2,2}_{loc}(U)$, see
\cite{GT, Theorem 8.8}. Hence by Sobolev embedding theorem,
$\Delta \Psi^I$ is in $L^p_{loc}(U)$ for $p=2n/(n-p)$, and
$\Psi^I\in W^{2,p}_{loc}(U)$, see \cite{GT, Lemma 9.16}. We can
then iterate by using the Sobolev embedding theorem to conclude
that the lemma is true.

To prove the claim, let $\Phi=\sum_I\Phi^I\sigma_I\in W^{1,2}_0(U)$. Then

$$
\nabla_{e_i}\Psi=\sum_{I}e_i(\Psi^I)
\sigma_I-\frac14\sum_{k,l,I}\Psi^I \omega_{kl}(e_i)  e_k\cdot
e_l\cdot \sigma_I
$$
and

$$
\nabla_{e_i}\Phi=\sum_{I}e_i(\Phi^I)
\sigma_I-\frac14\sum_{k,l,I}\Phi^I \omega_{kl}(e_i)  e_k\cdot
e_l\cdot \sigma_I
$$
where $\omega_{kl}$ are the connection forms with respect to the adapted frame $e_i$.
By Lemma 3.1

$$
\split \langle \nabla \Psi, \nabla \Phi \rangle=& \sum_I\langle
\nabla \Psi^I, \nabla \Phi^I \rangle
-\frac14\sum_{i=1}^{n-1}\sum_{ j,k,l,I,J}\langle \Psi^I \omega_{kl}(e_i)e_k \cdot e_l\cdot \sigma_I,e_i(\Phi^J)\sigma_J\rangle \\
&\quad-\frac14\sum_{i=1}^{n-1}\sum_{j,k,l,I,J } \langle
e_i(\Psi^I)\sigma_I, \Phi^J\omega_{kl}(e_i)e_k \cdot e_l \cdot
\sigma_J\rangle  +\sum_{I,J}a_{IJ}\Psi^I\bar\Phi^J
\endsplit
$$
where $a_{IJ}$ is a bounded function. Since $e_i(\omega_{kl}(e_i))$ is smooth up to boundary in $N\setminus \Omega$ and in $\Omega$, and for $1\le i\le n-1$ we can perform integration by parts to conclude that

$$
\sum_{i=1}^{n-1}\int_U\sum_{ j,k,l,I,J}\langle \Psi^I
\omega_{kl}(e_i)e_k \cdot e_l\cdot
\sigma_I,e_i(\Phi^J)\sigma_J\rangle
=\sum_{i=1}^{n-1}\sum_{IJ}\int_U\lf(e_i(\Psi^I)\ol\Phi^J
b_{iIJ}+\Psi^I\ol\Phi^J c_{iIJ}\ri),
$$
where  $b_{iIJ}$ and $c_{iIJ}$ are $L^\infty$ functions in $U$.
For simplicity, we set:

$$
\sum_{i=1}^{n-1}\int_U\sum_{ j,k,l,I,J}\langle e_i(\Psi^I)
\sigma_I,\Phi^J\omega_{kl}(e_i)e_k \cdot e_l\cdot \sigma_I\rangle
=\sum_{i=1}^{n-1}\sum_{IJ}\int_U\lf(e_i(\Psi^I)\ol\Phi^J
d_{iIJ}\ri),
$$
here, $d_{iIJ}$ are also $L^\infty$ functions in $U$.
By Lemma 3.2  and the fact that $\Psi$ is a
weak solution of $\D^2\Psi=0$, it is easy to see that the claim is true with 

$$\split &{A^J}_{Ii}=-\frac14(b_{iIJ}+d_{iIJ}),\\
&{B^J}_I =a_{IJ}+\sum_{i=1}^{n-1}c_{iIJ} +\frac{\Cal R}{4}\delta_{IJ} 
\endsplit
$$
where $\delta_{II}=1$ and $\delta_{IJ}=0$ if $I\neq J$. 
 This completes the proof of the lemma.

\enddemo

As a corollary, we have:

\proclaim{Corollary 3.1} Suppose $\Psi$ is $W^{1,2}$ weak solution
of

$$\D^2 \Psi =0 $$
in an open set $U$ in $N$. Then $\D\Psi \in W^{1,2}_{loc}(U)$.

\endproclaim

\demo{Proof} It is sufficient to consider the behavior of $\Phi$
near a point $x$ in $\p \Omega$. We choose an adapted orthonormal
frame near $x$ as before. With the notations as in the proof of
the previous lemma, we have

$$
\D\Psi = \sum_I\nabla \Psi^I \cdot \sigma_I -\frac14\sum_{i,j,k=1}^{n-1}\sum_{I}\omega_{kl}(e_i)e_i\cdot e_k \cdot
e_l \cdot \sigma_I - H e_n \cdot \Psi 
$$
where we have used Lemma 3.1. Here $H$ is the mean curvature of
the level surface $\rho$=constant. By Lemma 3.3,    first term in
the above equality is in $W^{1,2}_{loc}$. By Lemma 3.1, by the
assumption of the smoothness of $g$, we see for $1\le i, k,l\le
n-1$, $\omega_{kl}(e_i)$ is Lipschitz. By the assumption of the
mean curvature on $\partial \Omega$, we see $H$ is also Lipschitz
on $N$. The   corollary follows because being in $W^{1,2}_{loc}$
does not depend on the choice of orthonormal frame.

\enddemo

\proclaim{Lemma 3.4} Suppose $\sigma \in W^{1,2}_{loc}(N)$,
$\int_N \|\sigma\|^2  <\infty$ and $\D \sigma =0$. Then
$\sigma=0$.
\endproclaim
\demo{Proof} By the assumption and Lemma 3.2,  for any $\xi \in
W^{1,2}_0 (N)$, we have:
$$
0=\int_N \langle\ D \sigma, \D \xi\rangle =
\int_N \langle\nabla \sigma, \nabla \xi\rangle + \frac{1}{4}\int_N
\Cal R \langle \sigma, \xi\rangle
$$
Let $\xi=\eta^2 \sigma$, here $\eta$ is a cut-off function such
that for   $\rho$$>0$,

$$
\eta=\cases 1&\text{\ in $B_{\rho}(o)$}\\
0&\text{\ outside $B_{2\rho}(o)$},\endcases
$$
and 
$$
|\nabla \eta|\leq \frac{C}{\rho}.
$$
Then, we get:
$$
\int_{B_\rho} \|\nabla \sigma\|^2 \leq \frac{C}{\rho}\int_N \|\sigma \|^2 \to 0.
$$
Hence, $\sigma$ is parallel inside and outside $\Omega$.  Thus $\sigma =0$ outside $\Omega$ because  $\int_N||\sigma||^2<\infty.$ $\sigma=0$ inside $\Omega$ because $\sigma$ is continuous on $N$   by Lemma 3.3 and is parallel. This completes the proof of
the lemma. 
\enddemo

Now one can proceed as in the case that $g$ is smooth to prove the positive mass theorem.

\demo{Proof of Theorem 3.1} Let us first prove that $m_E\ge0$ for all end $E$. We assume that $N$ has only one end, the proof for the general case is similar. Let $\eta$ be a parallel spinor outside $\R^n$ with respect to the Euclidean metric. We may extend $\eta$ so that it is zero on a neighborhood of $\ol\Omega$. By the asymptotic conditions on $g$, we have $||\eta||$ is asymptotically constant,

$$
||\D\eta||(x)=O\lf(r^{1-n}(x) \ri),\tag3.3
$$
and
$$
||\D^2\eta||(x)=O\lf(r^{-n}(x)\ri).\tag3.4
$$
Here $r$ is the geodesic distance function with respect to $g$. Let $R>0$ be large enough,  then one can find spinor $\Psi_R\in W^{1,2}(B_o(R))$ where $o\in N$ is a fixed point, so that

$$
\D^2\Psi_R=0
$$
in the weak sense in $B_o(R)$ such that $\Phi_R=\eta$ on $\p B_o(R)$. This is equivalent to solve the following:

$$
\cases
 \D^2 \sigma_R &=-\D^2 \eta \text{\ \ in $B_o(R)$,}\\
\sigma_R|_{\p B_o(R)} &=0
\endcases\tag*
$$
One may use Lax-Milgram theorem to solve (*). Indeed, in the
Hilbert space consisting of all spinors in ${W_0}^{1,2}(B_o(R))$, define the sesqui-bilinear form:

$$
a(\Phi,\Psi)=\int_{B_o(R)}\langle \D\Phi,\D\Psi).
$$ 
Consider the linear functional

$$
F(\Psi)= -\int_{B_o(R)}\langle\D\Psi,\D\eta\rangle.
$$
It is easy to see that:

(i) $a$   is bounded, i.e. there is $C>0$ such that:
$$
|a(\Phi,\Psi)|\leq C|||\Phi|||\,|||\Psi|||.
$$

(ii)  It is positive by Lemma 3.2, the fact that $\Cal R\ge0$ and the Poincar\'e inequality. That is, there is $\delta$$>0$ such that :

$$
a(\Psi,\Psi) \geq \delta |||\Psi|||^2.
$$

(iii) $F$ is bounded.

Then by  Lax-Milgram theorem \cite{Yo, Sec.7, Chap.3}, we conclude that (*) has a solution.
 By Lemma 3.3, $\Psi_R$ is bounded. Hence $||\Psi_R||^2$ is in $W^{1,2}(B_o(R))$. Moreover, if $f\in C_0^\infty(B_o(R))$ with $f\ge0$, then

$$
\split \int_{B_o(R)}\langle \nabla ||\Psi_R||^2,\nabla f\rangle
&=\int_{B_o(R)}\lf(\langle \nabla \Psi_R,\nabla(f\Psi_R)\rangle+\langle \nabla (f\Psi_R),\nabla \Psi_R\rangle\ri)-2\int_{B_o(R)}f||\nabla \Psi_R||^2\\
&=-\frac12\int_{B_o(R)}f\lf(\Cal R||\Psi_R||^2+4||\nabla\Psi_R||^2\ri)\\
&\le0
\endsplit
$$
where we have used the Lichnerowicz formula in Lemma 3.2. Hence $||\Psi_R||^2$ is subharmonic in the weak sense. Since $\eta$ is uniformly bounded, we conclude that $\Psi_R$ are uniformly bounded by the maximum principle. Hence there is $R_i\to\infty$ such that $\Psi_i=\Psi_{R_i}$ converges in $W^{1,2}_{loc}(N)$ to $\Psi$ with $\D^2\Psi=0$ in the weak sense.

We claim that $\Psi$ is asymptotically close to $\eta$ in the following sense:

$$
||\Psi-\eta||(x)\le Cr^{2-n}(x)\log r(x)\tag3.5
$$
for some constant if $r(x)$ is large enough, and

$$
\int_{N}||\nabla
(\Psi-\eta)||^2+||\D(\Psi-\eta)||^2<\infty.\tag3.6
$$

To prove the claim, let us assume  $\ol\Omega \subset B_o (R_0)$ and that 
$R_i>R_0$ for all $i$. Then for any $i$, since $\Psi_i-\sigma=0$ on $\p B_o(R_i)$, we have

$$
\int_{B_o(R_i)}\langle \D\Psi_i,\D(\Psi_i-\sigma)\rangle=0.
$$
From this, it is easy to see that

$$
\int_{B_o(R_i)}||\D(\Psi_i-\eta)||^2\le \int_N||\D\sigma||^2.
$$
Moreover, by Lemma 3.2 and the fact that $\Cal R\ge0$, it is easy to see that

$$
\int_{B_o(R_i)}||\nabla (\Psi_i-\eta)||^2\le
\int_{B_o(R_i)}||\D(\Psi_i-\eta)||^2\le  \int_N||\D\sigma||^2.
$$
Using (3.3) and the fact that $\Psi_i$ converge  to $\Psi$ in $W^{1,2}_{loc}(N)$ weakly, we conclude that (3.6) is true.

On the other hand, since $\Psi_i$ are uniformly bounded, there is a constant $C_1$ such that

$$
||\Psi_i-\eta||\le C_1
$$
on $\p B_o(R_0)$. Let $u\ge0$ be a solution of $\Delta u\le -||\D^2\eta||$ outside $B_o(R_0)$ such that $u\ge C_1$ on $\p B_o(R_1)$ and such that $u(x)\le C_2r^{2-n}(x)\log r(x)$ for some constant $C_2>0$. Such $u$ can be found, see for example \cite{S}. Note that outside $B_o(R_0)$, $\Psi_i$ is smooth and by the usual Lichnerwicz formula,

$$
\Delta ||\Psi_i-\eta||\ge -||\D^2\eta||
$$
in the weak sense on $B_o(R_i)\setminus B_o(R_0)$. Hence

$$
||\Psi_i-\eta||(x)\le u(x)
$$
on  $B_o(R_i)\setminus B_o(R_0)$. Note that $\Psi_i$ converges pointwisely to $\Psi$ outside $B_o(R_0)$, hence (3.5) is true.

By Corollary 3.1, $\D\Psi\in W^{1,2}_{loc}(N)$. Apply Lemma 3.4 to $\sigma=\D\Psi$ and using (3.6), we conclude that $\D\Psi=0$.

Now choose $\eta$ to be nonzero constant spinor (with respect to
the Euclidean metric) near infinity normalized so that $||\eta|| $
is asymptotically 1 at infinity. Let $S(r)$ be the Euclidean
sphere with radius $r$ near infinity of $N$ and let $U(r)$ be the
interior of $S(r)$, then by Lemma 3.2,

$$
\int_{S(r)}\langle
\nu\cdot\D\Psi+\nabla_\nu\Psi,\Psi\rangle>=\int_{U(r)}||\nabla
\Psi||^2+\frac14\int_{U(r)}\Cal R||\Psi||^2.\tag3.7
$$
where $\nu$ is the outward normal of $S(r)$. On the other hand,

$$
\split \int_{S(r)}\langle
\nu\cdot\D\Psi+\nabla_\nu\Psi,\Psi\rangle>&=
\int_{S(r)}\langle  \D_T\Psi,\Psi\rangle>\\
&=\int_{S(r)}\langle  \D_T(\Psi-\eta),\Psi-\eta\rangle>+\int_{S(r)}\langle  \D_T(\Psi-\eta),\eta\rangle>\\
&\quad+\int_{S(r)}\langle  \D_T\eta,\Psi-\eta\rangle>+\int_{S(r)}\langle \D\eta,\eta\rangle\\
&=\int_{S(r)}\langle  \D_T(\Psi-\eta),\Psi-\eta\rangle>+\int_{S(r)}\langle  \Psi-\eta,\D_T\eta\rangle>\\
&\quad+\int_{S(r)}\langle
\D_T\eta,\Psi-\eta\rangle>+\int_{S(r)}\langle \D\eta,\eta\rangle
\endsplit\tag3.8
$$
where $\D_T=\sum_{i=1}^{n-1}e_i\nabla_{e_i}$, with $e_i$ to be
orthonormal and tangential to $ S(r)$. Here we have used the fact
that $D_T$ is self-adjoint on the boundary.

By (3.5) and (3.6), for each $r$ large enough,  we may
choose $r'$$\in$$(r, 2r)$, such that:

$$
\split\lf|\int_{\partial S(r')}\langle \D_T (\Psi-\eta), \Psi-\eta\rangle\ri|&\le   \int_{\partial B_r'}\|\langle \D_T (\Psi-\eta), \Psi-\eta\rangle\|\\
&=\frac 1r\int_{U(2r)\setminus U(r)}\|\langle \D_T (\Psi-\eta), \Psi-\eta\rangle\|\\
&\le \frac 1r\lf(\int_{U(2r)\setminus U(r)}\lf|\langle \D_T (\Psi-\eta)\ri|^2\ri)^\frac12\lf(\int_{U(2r)\setminus U(r)} ||\Psi-\eta||^2\ri)^\frac12\\
&\le  C r^{1-\frac n2}\log r.
\endsplit
$$
Thus, we see that we can find $r_i\to\infty$ such that

$$
\lim_{i\to\infty}\int_{\partial S(r_i)}\langle \D_T (\Psi-\eta),
\Psi-\eta\rangle=0.\tag3.9
$$
By (3.5) and (3.3), we also have

$$
\lim_{i\to\infty}\int_{S(r_i)}\langle
\Psi-\eta,\D_T\eta\rangle>+\int_{S(r_i)}\langle
\D_T\eta,\Psi-\eta\rangle>=0.\tag3.10
$$ 
Finally, by the argument in [B1, p.691--692], we can prove that

$$
\lim_{i\to\infty}\int_{S(r_i)}\langle \D\eta,\eta\rangle=C(n)m_E
$$
for some positive constant $C(n)$ depending only on $n$. Combining (3.7)-(3.10),
we conclude that $m_E\ge0$.

Suppose the mass of  some end $E_1$ is zero. Suppose $N$ has at
least two ends $E_1$, $E_2$. Then we may choose $\eta$ such that
$\eta$ is almost parallel but nonzero in $E_1$ and $\eta$ is zero
on $E_2$ and other  ends. By the above arguments, we conclude that
there is a spinor $\Psi$ which is asymptotically close to $\eta$
near infinity. Moreover, $\Psi$ is parallel, namely $\Psi$ is
smooth and parallel in the usual sense in the interior and
exterior of $\Omega$.  By Lemma 3.3, $\Psi$ is continuous. This is
impossible. Therefore $N$ has only one end. By choosing enough
linearly independent $\eta$ and by constructing continuous
parallel spinors from $\eta$, we can conclude as in \cite{B1} then
the curvature of $N$ is zero both inside and outside $\Omega$.
This completes the proof of the theorem.

\enddemo

\subheading{\bf \S4 Compact manifolds with boundary and with nonnegative scalar curvature}

In this section, we will use the results in sections 2 and 3 to study the boundary behaviors of a compact Riemannian manifold $(\Omega^n,g)$ of dimension $n$ with smooth boundary $\p\Omega$ and with  nonnegative scalar curvature. First we need the following lemma.

\proclaim{Lemma 4.1} Let $M$ be a smooth differentiable manifold and $\Omega$ be a  domain in $M$ with smooth boundary. Suppose $g$ is a Riemannian metric on $M$ satisfying the following:
\roster
\item"{(a)}"   $g_{M\setminus \ol\Omega}$ and $g_\Omega$ are smooth up to the boundary $\p\Omega$ and $g$ is Lipschitz near any point on $\p\Omega$.
\item"{(b)}" The sectional curvature of $M\setminus\ol\Omega$ and $\Omega$ are zero near  $\partial \Omega$.
\item"{(c)}"    $\partial \Omega$ has the same second
fundamental form with respect to $g_{M\setminus \ol\Omega}$ and $g_\Omega$ (with respect to the same normal direction).
\endroster  Then $g$ is $C^2$ in
a neighborhood of $\partial \Omega$.
\endproclaim
\demo{Proof} Let $p\in\p\Omega$ and let $\rho$ be the signed distance function from $\p\Omega$. Near $p$, the metric $g$ can be expressed in the form:

$$g=d\rho^2 + \sum_{i,j=1}^{n-1} g_{ij}d\theta_i d\theta_j,
$$
where 
$\sum_{i,j=1}^{n-1} g_{ij}(\rho,\theta)d\theta_i d\theta_j$ is the induced metric on the level sets of $\rho$ and $(\theta_1,\dots,\theta_{n-1},\rho)$ are the local coordinates. It is sufficient to prove that the components $g_{ij}$ are in $C^2$. Note that   partial derivatives of $g_{ij}$ of all order  with respect to $\theta$ exist. 

Let $\sum_{ij}^{n-1}
h_{ij}d\theta_i d\theta_j$ be the second fundamental form on the level surface $\rho$=constant with respect to the unit normal $\p/{\p \rho}$. Then for $\rho\neq0$, we have that

$$
\split
h_{ij}&=-\langle\nabla_{\frac{\partial}{\partial \theta_i}}\frac{\partial}{\partial \theta_j},\frac{\partial}{\partial \rho}\rangle\\
&=-\Gamma_{ij}^n\\
&=\frac12\frac{\p g_{ij}}{\p\rho}
\endsplit
$$
where $\Gamma_{ab}^c$ are the   Christoffel symbols and $\rho$ is considered to be the $n$-th  coordinate. Hence
$$\frac{\partial g_{ij}}{\partial \rho}=2h_{ij}.\tag4.1
$$
By the assumption that $h_{ij}$ agrees on $\p\Omega$, we see that $\p g_{ij}/\p\rho$ is continuous up to $\p\Omega$. Hence $g_{ij}$ is $C^1$ near $p$.

By (4.1), for $\rho\neq0$
$$\split
\frac{\partial^2 g_{ik}}{\partial\rho \partial\theta_j}&=\frac{\partial^2 g_{ik}}{\partial\theta_j \partial\rho} \\
&=2\frac{\partial{h_{ik}}}{\partial\theta_j}.
\endsplit
$$
Since $h_{ij}$ agrees on $\p\Omega$ and $g$ is smooth up to the boundary when restricted on $\Omega$ or on $M\setminus\ol\Omega$, we conclude that ${\partial^2 g_{ik}}/{\partial\rho \partial\theta_j}$ is continuous near $p$.  

Next, we want to show $\frac{\partial^2 g_{ik}}{\partial \rho^2}$
is also continuous. For $\rho\neq0$, using the fact that the sectional curvature is zero near $p$, we have
$$
\split
\frac{\p h_{ij}}{\p\rho}&=-\frac{\p}{\p\rho}\langle\nabla_{\frac{\partial}{\partial \theta_i}}\frac{\partial}{\partial \theta_j},\frac{\partial}{\partial \rho}\rangle\\
&=-\langle\nabla_{\frac{\p}{\p\rho}}\nabla_{\frac{\p}{\p \theta_i}}\frac{\p}{\p \theta_j},\frac{\p}{\p  \rho}\rangle\\
&=-\langle\nabla_{\frac{\p}{\p \theta_i}}\nabla_{\frac{\p}{\p\rho}}\frac{\p}{\p \theta_j},\frac{\p}{\p  \rho}\rangle\\
&=-\frac{\p}{\p\theta_i}\langle\nabla_{\frac{\p}{\p\rho}}\frac{\p}{\p \theta_j},\frac{\p}{\p  \rho}\rangle+\langle\nabla_{\frac{\p}{\p\rho}}\frac{\p}{\p \theta_j},\nabla_{\frac{\p}{\p \theta_i}}\frac{\p}{\p  \rho}\rangle\\
&=\langle\nabla_\frac{\p}{\p \theta_j}{\frac{\p}{\p\rho}},\nabla_{\frac{\p}{\p \theta_i}}\frac{\p}{\p  \rho}\rangle\\
&=g^{ks}g^{ls}h_{ik}h_{is}
\endsplit
$$
where $(g^{ij})$ is the inverse of $(g_{ij})$. Combining this with (4.1), we see that $\frac{\partial^2 g_{ij}}{\partial \rho^2}$ is  continuous. Hence $g$ is in $C^2$.
\enddemo
 \proclaim{Remark}  We can prove $g$ is actually $C^{\infty}$ by the
same argument.
\endproclaim

 Let $(\Omega^n,g)$ be a Riemannian manifold of dimension $n$ with compact closure and with smooth boundary. Let us first consider the case that $n>3$. We assume the following:

\roster

\item"{(i)}" $\p\Omega$ has finitely many components $\Sigma_i$, $1\le i\le k$.

\item"{(ii)}" The mean curvature $H$ of $\Sigma_i$ with respect to the outward normal is positive.

\item"{(iii)}"  There is an isometric embedding $\iota_i:  \Sigma_i\to\R^n$ such that $\Sigma_i$ is a strictly convex closed hypersurface in $\R^n$. Here we identify $(\Sigma_i,g)$ with its image with the metric induced by the Euclidean metric in $\R^n$.\item"{(iv)}" $\Omega$ is spin.

\endroster

We may extend $\Omega$ across $\p\Omega$ to a smooth manifold $\wt\Omega$ which contains $\ol\Omega$. By the embedding of $\Sigma_i$, we can define a diffeomorphism from a neighborhood of $\Sigma_i$ in $\wt\Omega$ to a neighborhood of $\Sigma_i$ in $\R^n$ by mapping the set with distance $r$ from $\Sigma_i$ in $\wt\Omega$ to the set with distance $r$ from $\Sigma_i$ in $\R^n$, so that the part near $\Sigma_i$ which is outside of $\Omega$ in $\wt\Omega$ will be mapped into an open set which is outside of $\Sigma_i$ in $\R^n$.

\proclaim{Theorem 4.1} Let $(\Omega^n,g)$ be a compact manifold
with smooth boundary and with nonnegative scalar curvature and $n>3$.
Suppose $\Omega$ satisfies conditions (i)-(iv). Then for each boundary
component $\Sigma_i$,

$$
\int_{\Sigma_i}Hd\sigma\le \int_{\Sigma_i}H^{(i)}_0d\sigma\tag4.2
$$
where $H^{(i)}_0$ is the mean curvature of $\Sigma_i$ in $\R^n$
with respect to the outward normal. Moreover, if equality holds in
(4.2) for some $i$, then $\p\Omega$ has only one boundary
component and $\Omega$ is a domain in $\R^n$.

\endproclaim

To prove the theorem, let us fix some notations. For each $i$, we may suppose $\Sigma_i$ is a strictly convex hypersurface in $\R^n$. For simplicity, let us denote $\Sigma_i$ by $\Sigma_0$ and $H_0^{(i)}$ by $H_0$. In the setting as in section 2, let $u$ be the solution of (2.1) with initial data $u(x,0)=H_0(x)/H(x)$ which is positive by (ii) and (iii).

\proclaim{Lemma 4.2} The function

$$
m(r)=\int_{\Sigma_r}H_0(1-u^{-1})d\sigma_r
$$
is nonincreasing in $r$, where $H_0$ is the mean curvature of $\Sigma_r$ in $\R^n$.

\endproclaim

\demo{Proof} Let $h_{ij}^0$ be the second fundamental form of $\Sigma_r$ with respect to the Euclidean metric. Then by the Gauss equations, it is easy to see that

$$
\frac{\p H_0}{\p r}=-\sum_{i,j=1}^{n-1}\lf(h_{ij}^0\ri)^2.
$$

Since  $u$ satisfies (2.1), we get:

$$
\split
\frac{d}{dr}&\int_{\Sigma_r}H_0 (1-u^{-1})d\sigma_r\\
&=-\int_{\Sigma_r}\sum_{i,j=1}^{n-1}
\lf(h^0_{ij}\ri)^2
(1-u^{-1})d\sigma_r+\int_{\Sigma_r}u^{-2}H_0\frac{\partial u}{\partial r}  +\int_{\Sigma_r}H_0^2 (1-u^{-1})d\sigma_r\\
&=\int_{\Sigma_r}\lf[\lf(H_0^2-\sum_{i,j}\lf(h^0_{ij}\ri)^{2}\ri)(1-u^{-1})
+\Delta_r u+\frac12(u^{-1}-u)\Cal R^r\ri] d\sigma_r \\
&=-\frac12\int_{\Sigma_r}\Cal R^r u^{-1}(1-u)^2\\
&\leq 0.
\endsplit
$$
where $\Delta_r$ is the Laplacian on $\Sigma_r$ and  $\Cal R^r$ is the scalar curvature of $\Sigma_r$ with respect to the induced metric in $\R^n$ (which is the same as the metric induced by $g$). Here  we have used the fact

$$H_0^2-\sum_{i,j}(h^0_{ij})^{2}= \Cal R^r.
$$
Thus, we see that $m(r)$ is nonincreasing.

\enddemo

We are ready to prove the theorem.

\demo{Proof of Theorem 4.1} Using the above method, we can attach each boundary component
$\Sigma_i$ to the exterior of a convex hypersurface in $\R^n$, which is denoted by $E_i$.
 On each $E_i$, we construct the metric given by $g_i=u_i^2 dr^2 +g_r$ with initial data
  $H_0^{(i)}/H$ as in Theorem 2.1. Denote the resulting manifold by $N$. Let $g_N$ be the
   metric on $N$ defined by $g_N=g$ in $\Omega$ and $g_N=g_i$ on each $E_i$. Since $\Omega$
   is spin, $N$ is spin. By Theorem 2.1 and (1.6), $N$ satisfy the assumptions in Theorem 3.1.
   Hence the mass $m_{E_i}$ is nonnegative for each $i$. By Theorem 2.1(c) and Lemma 4.1,
    we conclude that (4.2) is true for all $i$.

Suppose equality holds in (4.2) for some $i$, then the mass of
$E_i$ must be zero. Hence $N$ has only one end and $N$ is flat by
Theorem 3.1. Therefore $\p\Omega$ has only one component and
$\Omega$ is flat. By (1.7), we have $u\equiv1$.  On the other hand, we note that $n>3$  and the
boundary is strictly convex in $\R^n$. By  the proof in  \cite{E,\S60}, we know that the second fundamental forms of
the boundary of $\Omega$  with respect to $g$ and the Euclidean metric (in the  same normal direction) are equal.
By Lemma 4.1, we see the metric on $N$ is actually $C^2$. Since $u\equiv1$, $N$ is the Euclidean space outside a compact set. By volume comparison theorem,  $N$ is
isometrically to $\R^n$ which implies that $\Omega$ is a domain in
$\R^n$. This completes the proof of the theorem.
\enddemo

In case $n=3$ then  condition (iv) mentioned above is automatically satisfied. Also, by a well-known result, see  \cite{N} for example, condition (iii) is equivalent to the condition that $\Sigma_i$ has positive Gaussian curvature. It is also well-known that  the embedding is unique up to an isometry in $\R^3$. Hence in this case we have:

\proclaim{Theorem 4.2} $(\Omega^3,g)$ be a Riemannian manifold of dimension $3$ with compact closure  with smooth boundary and with
 nonnegative scalar curvature. Suppose $\Omega$ satisfies conditions (i)-(ii).
Moreover, suppose each boundary component of $\Sigma_i$ has positive Gaussian curvature.
 For each boundary component $\Sigma_i$, we have

$$
\int_{\Sigma_i}Hd\sigma\le \int_{\Sigma_i}H^{(i)}_0d\sigma\tag4.3
$$
where $H^{(i)}_0$ is the mean curvature of $\Sigma_i$ with respect to the outward normal
when it is isometrically embedded in $\R^3$. Moreover, if equality holds in (4.3) for some $i$,
 then $\p\Omega$ has only one   component and $\Omega$ is a domain in $\R^3$.

\endproclaim

\demo{Proof}  The proof of (4.3) is the same as before. Using the same notations as in the proof of Theorem 4.1,  by the same argument as in the proof of Theorem 4.1, if
equality holds in (4.3) for some end, then $\p\Omega$ has only one component and $N$ is flat with $u\equiv1$. It remains to prove that
$\Omega$ is actually a domain in $\R^3$. 

Since $u\equiv1$,  the mean curvatures of $\p \Omega$ with respect to $g$ and the Euclidean metric  are  equal. Since $\p\Omega$ is now a strictly convex surface in $\R^3$, it is easy to see that   the proof in \cite{K, Theorem 6.2.8} can be carried over and we can conclude that  
the second fundamental forms of $\p\Omega$ with respect to
$g$ and the Euclidean metric  are equal. Hence we can conclude from Lemma 4.1 as before that   $\Omega$
is a domain in $\R^3$. This finishes the proof of the theorem.
\enddemo

By the result of Weyl \cite{We},  the boundary component $\Sigma_i$ of $\Omega$ in Theorem 4.2 satisfies

$$
4K\le \lf(H_0^{(i)}\ri)^2\le \sup_{\Sigma_i}\lf(4K-K^{-1}\Delta
K\ri)$$
where $\Delta$ is the Laplacian of $\Sigma_i$ with metric induced by $g$, $K$ is the Gaussian curvature of $\Sigma_i$ and $H_0^{(i)}$ is the mean curvature of $\Sigma_i$ when it is embedded in $\R^3$. Hence we have the following corollary.

\proclaim{Corollary 3.1} Let $(\Omega^3,g)$ be a compact manifold of dimension 3 with boundary and with nonnegative scalar curvature. Suppose $\Omega$ satisfies conditions (i)-(ii). Moreover, suppose each boundary component of $\Sigma_i$ has positive Gaussian curvature $K$. Then

$$
\frac{1}{Area(\Sigma_i)} \int_{\Sigma_i}Hd\sigma\le
\lf[\sup_{\Sigma_i}\lf(4K-K^{-1}\Delta K\ri)\ri]^\frac12. $$
Moreover, if equality holds for some $\Sigma_i$, then $\p\Omega$ has only one component and $\Omega$ is a domain in $\R^3$.

\endproclaim
\subheading{\bf \S5 An equivalent statement of the positive mass theorem}

In the previous section, we obtain Theorem 4.2 from the positive mass theorem: Theorem 3.1. In this section, we want to show that one can obtain the first part of the positive mass theorem by assuming that Theorem 4.2 is true. More precisely, let $(N,g)$ be a complete noncompact manifold with finitely many ends with the following properties:
\roster
\item"{(i)}" $N$ has nonnegative scalar curvature $\Cal R$  which is in $L^1(N)$. 
\item"{(ii)}" Each end $E$ is diffeomorphic to the exterior of some compact set in $\R^3$, so that the metric $g=\sum_{i,j=1}^3g_{ij}dx_idx_j$ is   asymptotically flat in the sense that
$$
g_{ij}=\delta_{ij}+\tilde b_{ij}
$$
such that 
$$|\tilde b_{ij}|+r|\nabla_0\tilde b_{ij}|+r^2|\nabla_0\nabla_0 \tilde b_{ij}|=O(r^{-1}).$$ 
where   $r$ is the Euclidean distance from the origin and $\nabla_0$ is the derivatives with respect to the Euclidean metric.
\endroster 
Then we have:

\proclaim{Theorem 5.1} Suppose (4.2) is true for any compact Riemannian three manifold $\Omega$ with boundary satisfying the assumptions in  Theorem 4.2.  Let $(N,g)$ be as above, then the ADM mass of each end of $N$ is nonnegative. 
\endproclaim

\demo{Proof} By the result of \cite{SY5}, it is sufficient to prove the theorem   under the stronger assumption  that at each end $E$,
$$
g_{ij}=\lf(1+\frac{2m_E}r\ri)\delta_{ij}+b_{ij},\tag5.1
$$

$$|b_{ij}|+r|\nabla_0 b_{ij}|+r^2|\nabla_0\nabla_0 b_{ij}|+r^3|\nabla_0\nabla_0\nabla_0 b_{ij}|+r^4|\nabla_0\nabla_0\nabla_0\nabla_0 b_{ij}|=O(r^{-2})$$ 
where $m_E$ is a constant, and $r$ is the Euclidean distance from the origin. Namely, it is sufficient to prove that $m_E\ge0$ under the additional condition (5.1).

 In fact, by \cite{SY5}, given any $\e>0$, one can construct a new metric $\wt g$ on $N$ with zero scalar curvature such that near infinity at each end $E$, the metric $\wt g$ is of the form:
$$
\wt g_{ij}=\varphi^4 \delta_{ij},
$$
where $\varphi$ satisfies
$$
\varphi=1+\frac{\wt m_E}r+h.
$$
Here $\wt m_E$ is a constant and $|h|=O(r^{-2})$. Moreover, $\wt m_E\le m_E+\e$. Since the scalar curvature of $\wt g$ is zero, $\varphi$ and hence $h$ is harmonic outside a compact set of the Euclidean space. By the gradient estimates of harmonic functions on Euclidean space, we conclude that $\wt g$ satisfies (5.1).

We assume that $N$ has only one end, and denote $m_E$ by $m$. The general case can be proved similarly. The proof of Theorem 5.1 are divided into several steps. 

{\sl Step 1:} We want to compute   the Gaussian curvature $K$ of $\Bbb S(r)$ with respect to the metric $g$. Here $\Bbb S(r)$ is the Euclidean sphere of radius $r$. 

Let $\By=\By(\z_1,\z_2)$ be local parametrization of the standard unit sphere $\Bbb S=\Bbb S(1)$, here $\By=(y_1,y_2,y_3)$. Then local parametrization  for the standard $\Bbb S(r)$ is given by
$$
\bx=r\By.
$$
Here and below, $i,j...$ are from 1 to 3, and $\a,\b...$ are from 1 to 2. Now
$$
\frac{\p x_i}{\p r}=y_i\tag5.2
$$
and
$$
\frac{\p x_i}{\p \z_\a}=r\frac{\p y_i}{\p \z_\b}.\tag5.3
$$

 Let us first compute the metric on $\Bbb S(r)$. Since
$$
\frac{\p}{\p \z_\a}=\frac{\p x_i}{\p \z_a}\frac{\p}{\p x_i},
$$
we have 
$$
\split
\tau_{\ab}&=g\lf(\frac{\p}{\p\z_\a},\frac{\p}{\p\z_\b}\ri)\\
&=g_{ij}\frac{\p x_i}{\p \z_\a}\frac{\p x_j}{\p \z_\b}\\
&= r^2\lf[\lf(1+\frac{2m}r\ri)\delta_{ij}+b_{ij}\ri]\frac{\p y_i}{\p \z_\a}\frac{\p y_j}{\p \z_\b}\\
&=r^2\lf[\lf(1+\frac{2m}r\ri)a_\ab+b_{ij}\frac{\p y_i}{\p \z_\a}\frac{\p y_j}{\p \z_\b}\ri] 
\endsplit\tag5.4
$$
where $a_\ab$ is the standard metric on $\Bbb S(1)$ in the coordinates $(\z_1,\z_2)$. Hence 
$$
\split
\tau&=\tau_{11}\tau_{22}-\tau_{12}^2\\
&=ar^4\lf(1+\frac{2m}r\ri)^2(1+f)
\endsplit\tag5.5
$$
where   $a=\det(a_{\ab})$ and  $f$ is a smooth function which satisfies
$$
|f|+|\p f|+|\p^2 f|+|\p^3 f|+|\p^4 f|=O\lf(r^{-2}\ri).\tag5.6
$$
$\p$ denotes the partial derivatives with respect to $\zeta_\a$. Moreover, $|\frac{\p f}{\p r}|=O(r^{-3})$.  Here and below $f$ always denotes a smooth function, but the meaning of $f$ may vary from line to line. The function $f$ in (5.5) satisfies (5.6) because of the assumptions on $b_{ij}$ and (5.3). The inverse of $(\tau_\ab)$ is given by
$$
\tau^{\ab}=r^{-2}\lf(1+\frac{2m}r\ri)^{-1}\lf(a^\ab+f\ri)\tag5.7
$$
where $(a^\ab)=(a_\ab)^{-1}$ and $f$ also satisfies (5.6).  Let $\Gamma_\ab^\gamma$ and $\wt\Gamma_\ab^\gamma$ be the Christoffel symbols of $\Bbb S(r)$ with induced metric $\tau_\ab$ and that of the standard unit sphere $\Bbb S$ in the coordinates $\z$. Then
$$
 \Gamma_\ab^\gamma=\wt\Gamma_\ab^\gamma+f\tag5.7
$$
with
$$
|f|+|\p f|+|\p^2 f|+|\p^3 f|=O\lf(r^{-2}\ri),
$$
and
$$
\frac{\p \Gamma_\ab^\gamma}{\p\z_\delta}=\frac{\p\wt \Gamma_\ab^\gamma}{\p\z_\delta}+f \tag5.8
$$
with 
$$
|f|+|\p f|+|\p^2 f|=O\lf(r^{-2}\ri).
$$
Hence the Gaussian curvature $K$ of $\Bbb S(r)$ with metric induced by $g$ is 
$$
\split
K&=-\tau_{11}^{-1}\lf[\lf(\Gamma_{12}^2\ri)_1-\lf(\Gamma_{11}^2\ri)_2+\Gamma_{12}^1\Gamma_{11}^2+\Gamma_{12}^2\Gamma_{12}^2-\Gamma_{11}^2\Gamma_{22}^2-\Gamma_{11}^1\Gamma_{12}^2\ri)\\
&=-\tau_{11}^{-1}\lf(-a_{11}+f\ri)\\
&=r^{-2}\lf(1+\frac{2m}r\ri)^{-1}\lf(1+f\ri) 
\endsplit\tag5.9
$$
with
$$
|f|+|\p f|+|\p^2 f|=O\lf(r^{-2}\ri).
$$

{\sl Step 2:} We want to show that the mean curvature $H$ on $\Bbb S(r)$ is positive for $r$ large enough. Moreover, we want to compute the integral of $H$ over $\Bbb S(r)$.

Let $H$ be the mean curvature of $\Bbb S(r)$ with respect to the outward normal $\bn$ in the metric $g$. Let $A(r)$ be the area of $\Bbb S(r)$. 
Then by the first variational formula and (5.5)

$$
\split
\int_{\Bbb S(r)} \langle \frac{\p}{\p r},\bn\rangle H d\sigma_r&=A'(r)\\
&=\lf(2r+2m\ri)\int_{\Bbb S(1)}(1+f)^\frac12 d\sigma_0+O(r^{-1})\\
&=8\pi(r+m)+O(r^{-1})
\endsplit$$
where $d\sigma_r$ is the volume form of $\Bbb S(r)$ with metric induced by $g$ and $d\sigma_0$ is the volume form of the standard unit sphere. 

Note that
$$
\frac{\p}{\p r}=\frac{\p x_i}{\p r}\frac{\p}{\p x_i}=\frac{x_i}{r}\frac{\p}{\p x_i} 
$$
and the gradient of $r$ with respect  to $g$ is 
$$\nabla r=g^{ij}\frac{\p r}{\p x_i}\frac{\p}{\p x_j}=g^{ij}\frac{x_i}r\frac{\p}{\p x_j}
$$
one obtains 
$$
\langle \frac{\p}{\p r},\nabla r\rangle=1\tag5.10
$$
and
$$
|\nabla r|^2=g^{ij}\frac{\p r}{\p x_i}\frac{\p r}{\p x_j}=1-\frac{2m}r+h\tag5.11
$$
where $h=O(r^{-2})$. 
Since $\Bbb S(r)$ is the level surface of the function $r$,  $\bn =|\nabla r|^{-1}\nabla r$. We have
$$
\int_{\Bbb S(r)} |\nabla r|^{-1} H d\sigma_r=\int_{\Bbb S(r)} \langle \frac{\p}{\p r},\bn\rangle H d\sigma_r=8\pi(r+m)+O(r^{-1})\tag5.12
$$

For any point $x$ on $\Bbb S(r)$, choose an orthonormal frame $e_i$ with respect to $g$ such that $e_1$, $e_2$ are tangential and $e_3$ is the unit outward normal. Moreover, assume the second fundamental form $h_{\ab}$ is diagonalized at $x$. By the Gauss equation, 
$$
h_{11}h_{22}=K-R_{1212}
$$
where $R_{ijij}$ is the curvature tensor of $N$. By the asymptotic behavior of $g$, we have $|R_{ijij}|=O(r^{-3})$. By (5.9), we conclude that
$$
h_{11}h_{22}>0,
$$
if $r$ is large enough and so $h_{11}$ and $h_{22}$ are of the same sign. Hence $H$ is either everywhere positive or everywhere negative when $r$ is large. By (5.12), we must have $H>0$. 

Since $H>0$ for $r$ large enough, we may apply mean value theorem and use (5.11), we have
$$
\int_{\Bbb S(r)} H d\sigma_r=8\pi r+O(r^{-1}).\tag5.13
$$

Since the Gaussian curvature of $\Bbb S(r)$ is positive for $r$ large enough,   $(\Bbb S(r),g)$ can be embedded isometrically in $\Bbb R^3$. Let $H_0$ be its mean curvature in $\R^3$.

{\sl Step 3}: We want to estimate $H_0$. Note that $H_0^2\ge 4K$, and by \cite{We}, we have
$$
4K\le H_0^2\le \sup_{\Bbb S(r)}\lf(4K-K^{-1}\Delta K\ri),\tag5.14
$$
where $\Delta$ is the Laplacian of $\Bbb S(r)$ with metric induced by $g$. $$
\split
K^{-1}\Delta K&=K^{-1}\lf[\tau^\ab\frac{\p^2 K}{\p\z_\a\p\z_\b}+\frac1{\sqrt\tau}\frac{\p}{\p\z_\a}\lf(\sqrt\tau \tau^\ab\ri)\frac{\p K}{\p \z_\b}\ri]\\
&=O(r^{-4}) 
\endsplit\tag5.15
$$
because by (5.9)
$$
\frac{\p^2K}{\p\z_\a\p\z_\b}=r^{-2}\lf(1+\frac{2m}r\ri)^{-1}\frac{\p^2f}{\p\z_\a\p\z_\b}=O(r^{-4}).
$$
Combining (5.9), (5.14) and (5.15) we have
$$
H_0=\frac 2r-\frac{2m}{r^2}+O(r^{-3}).
$$
Hence by (5.5)
$$
\split
\int_{S(r)}H_0d\sigma_r&=\lf(\frac 2r-\frac{2m}{r^2}+O(r^{-3})\ri)\cdot r^2
\lf(1+\frac{2m}r\ri)\lf(4\pi+O\lf(r^{-2}\ri)\ri)\\
&=8\pi(r+m)+O(r^{-1}).
\endsplit\tag5.16
$$

{\sl Step 4}: We can now conclude the proof of the theorem. Let $\Omega_r$ be the domain in $N$ so that $\p\Omega_r$ is $\Bbb S(r)$. Then $\Omega_r$ has nonnegative scalar curvature, $\p\Omega_r$ has positive mean curvature and positive Gaussian curvature. By assumptions, (4.2) is true for $\p\Omega_r$. Combining with (5.13) and (5.16) we have
$$
\split
0&\le \int_{\Bbb S(r)}\lf(H_0-H\ri)d\sigma_r\\
&=8\pi m+O(r^{-1}).
\endsplit
$$
Let $r\to\infty$, we have $m\ge0$.
\enddemo

\end